%% file: SerkovDA-DGA-2023.tex
\documentclass[sn-mathphys,Numbered,pdflatex]{sn-jnl}

\usepackage{graphicx}%
\usepackage{multirow}%
\usepackage{amsmath,amssymb,amsfonts}%
\usepackage{amsthm}%
\usepackage{mathrsfs}%
\usepackage[title]{appendix}%
\usepackage{xcolor}%
\usepackage{textcomp}%
\usepackage{manyfoot}%
\usepackage{booktabs}%
\usepackage{algorithm}%
\usepackage{algorithmicx}%
\usepackage{algpseudocode}%
\usepackage{listings}%
\usepackage{euscript}
\usepackage{upgreek}
\usepackage{stmaryrd} 

\theoremstyle{thmstyleone}%
\newtheorem{theorem}{Theorem}
\theoremstyle{thmstyletwo}%
\newtheorem{remark}{Remark}%

\theoremstyle{thmstylethree}%

\newtheorem{lem}{Lemma}
\newtheorem{cor}{Corollary}

\renewcommand{\proof}{\par\mbox{P r o o f.}\ \ }
\newcommand{\nint}[2]{{{#1}..{#2}}}
\renewcommand{\ge}{\geqslant}
\renewcommand{\le}{\leqslant}
\newcommand{\myles}{\mathrel\prec}
\newcommand{\myor}{\ensuremath{\lor}}
\newcommand{\myand}{\ensuremath{\mathrel{\&}}}

\newcommand{\stmb}[2]{{{\bf\overline{s}}(#1,#2)}}
\newcommand{\stm}[2]{{{\bf s}(#1,#2)}}

\newcommand\argmin{\operatornamewithlimits{\mathrm{argmin}}}
\newcommand\argmax{\operatornamewithlimits{\mathrm{argmax}}}

\newcommand{\beq}{\begin{equation}}
\newcommand{\beqnt}{\begin{equation}\notag}
\newcommand{\eeq}{\end{equation}}
\newcommand{\mydef}{\mathrel{\triangleq}}
\newcommand{\icP}[1]{\EuScript{P}({#1})}
\newcommand{\icPp}[1]{\EuScript{P^\prime}({#1})}
\newcommand{\myll}{\ensuremath{\forall}}
\newcommand{\res}[2]{\ensuremath{(#1|#2)}}
\newcommand{\myimp}{\ensuremath{\mathrel{\Rightarrow}}}
\newcommand{\myeqv}{\ensuremath{\mathrel{\Leftrightarrow}}}
\newcommand{\myemp}{\ensuremath{\varnothing}}
\newcommand{\myxst}{\ensuremath{\exists}}
\newcommand{\myLe}{\mathrel{\sqsubseteq}}
\newcommand{\mdef}{\ensuremath{\mathrel{\text{def}}}}
\newcommand{\mydefeq}{\stackrel{\mdef}{\myeqv}}

\newcommand{\RA}{\ensuremath{\mathbb R}}
\newcommand{\NA}{\ensuremath{\mathbb N}}

\newcommand{\naY}{\ensuremath{\mathrm{Y}}}
\newcommand{\naX}{\ensuremath{\mathrm{X}}}
\newcommand{\naOm}{\ensuremath{\Omega}}
\newcommand{\naZ}{\ensuremath{\mathrm{Z}}}
\newcommand{\naZo}[2]{\naZ\setbox0=\hbox{\text{$#1$}}\ifdim\wd0=0pt\else\res{#1}{#2}\fi}
\newcommand{\naT}{\ensuremath{\mathrm{T}}}
\newcommand{\naTc}{\ensuremath{\mathcal{T}}}
\newcommand{\naGE}[3]{\mathrm{(GE)}\setbox0=\hbox{\text{$#1$}}\ifdim\wd0=0pt\else[{#1};{#2}\,|\,{#3}]\fi}
\newcommand{\naOmo}[2]{\naOm\setbox0=\hbox{\text{$#1$}}\ifdim\wd0=0pt\else\res{#1}{#2}\fi}
\newcommand{\nagH}[3]{\gamma_{#1}\setbox0=\hbox{\text{$#2$}}\ifdim\wd0=0pt\else[{#2}]\setbox0=\hbox{\text{$#3$}}\ifdim\wd0=0pt\else({#3})\fi\fi}
\newcommand{\nagi}[3]{\gamma^\infty_{#1}\setbox0=\hbox{\text{$#2$}}\ifdim\wd0=0pt\else[{#2}]\setbox0=\hbox{\text{$#3$}}\ifdim\wd0=0pt\else({#3})\fi\fi}
\newcommand{\nagk}[3]{\gamma^k_{#1}\setbox0=\hbox{\text{$#2$}}\ifdim\wd0=0pt\else[{#2}]\setbox0=\hbox{\text{$#3$}}\ifdim\wd0=0pt\else({#3})\fi\fi}
\newcommand{\naG}[2]{\Gamma\setbox0=\hbox{\text{$#1$}}\ifdim\wd0=0pt\else[{#1}]\setbox0=\hbox{\text{$#2$}}\ifdim\wd0=0pt\else({#2})\fi\fi}
\newcommand{\naGi}[2]{\Gamma^\infty\setbox0=\hbox{\text{$#1$}}\ifdim\wd0=0pt\else[{#1}]\setbox0=\hbox{\text{$#2$}}\ifdim\wd0=0pt\else({#2})\fi\fi}
\newcommand{\naGk}[2]{\Gamma^k\setbox0=\hbox{\text{$#1$}}\ifdim\wd0=0pt\else[{#1}]\setbox0=\hbox{\text{$#2$}}\ifdim\wd0=0pt\else({#2})\fi\fi}
\newcommand{\naD}[2]{{\mathbf D}\setbox0=\hbox{\text{$#1$}}\ifdim\wd0=0pt\else[{#1}]\setbox0=\hbox{\text{$#2$}}\ifdim\wd0=0pt\else({#2})\fi\fi}
\newcommand{\naM}[1]{\mathfrak{N}\setbox0=\hbox{\text{$#1$}}\ifdim\wd0=0pt\else[{#1}]\fi}
\newcommand{\naMo}[2]{\mathfrak{N}_0\setbox0=\hbox{\text{$#1$}}\ifdim\wd0=0pt\else[{#1};{#2}]\fi}
\newcommand{\naN}{\AnaN{\naTc}{}}
\newcommand{\naNo}[1]{\naN\setbox0=\hbox{\text{$#1$}}\ifdim\wd0=0pt\else[{#1}]\fi}
\newcommand{\AnaN}[2]{\mathbf{N}_{(#1)}\setbox0=\hbox{\text{$#2$}}\ifdim\wd0=0pt\else\left[{#2}\right]\fi}
\newcommand{\AnaNO}[2]{\mathbf{N}^0_{(#1)}\setbox0=\hbox{\text{$#2$}}\ifdim\wd0=0pt\else\left[{#2}\right]\fi}
\newcommand{\nado}[1]{(\mathrm{DOM})\setbox0=\hbox{\text{$#1$}}\ifdim\wd0=0pt\else\left[{#1}\right]\fi}
\newcommand{\naNO}[1]{\naN^0\setbox0=\hbox{\text{$#1$}}\ifdim\wd0=0pt\else[{#1}]\fi}
\newcommand{\sres}[2]{\llparenthesis#1|#2\rrparenthesis}
\newcommand{\Ana}[2]{\ensuremath{\langle{#1}\rangle_{#2}}}
\newcommand{\fix}[1]{\ensuremath{{\mathbf{Fix}({#1})}}}
\newcommand{\UB}{{\mathbf U}}
\newcommand{\VB}{{\mathbf V}}
\newcommand{\UBB}{{\bar{\mathbf U}}}
\newcommand{\VBB}{{\bar{\mathbf V}}}
\newcommand{\dup}{\mathrm{d}}
\newcommand{\fref}[1]{{\rm(\ref{#1})}}   

\begin{document}

\title[On a construction of a partially non-anticipative multiselector]{On a construction of a partially non-anticipative multiselector and its applications to dynamic optimization problems}


\author*[1,2]{\fnm{Dmitrii} \sur{Serkov}}\email{d.a.serkov@gmail.com, serkov@imm.uran.ru}

\affil*[1]{\orgname{Krasovskii Institute of Mathematics and Mechanics}, \orgaddress{\street{S. Kovalevskaya Str., 16}, \city{Ekaterinburg}, \postcode{620990}, \country{RF}}}

\affil[2]{\orgname{Ural Federal University}, \orgaddress{\street{Mira Str., 32}, \city{Ekaterinburg}, \postcode{10587},\country{RF}}}


\abstract{Let sets of functions $\naZ$ and $\Omega$ on the time interval $\naT$ be given, let there also be a multifunction (m/f) $\alpha$ acting from $\Omega$ to $\naZ$ and a finite set $\Delta$ of moments from $\naT$.
The work deals with two questions:
the first one is the connection between the possibility of stepwise construction (specified by $\Delta$) of a value $z$ of $\alpha(\omega)$ for an unknown step-by-step implemented argument $\omega\in\Omega$ and the existence of a multiselector $\beta$ of the m/f $\alpha$ with a non-anticipatory property of special kind defined by $\Delta$;
and the second question is how to build the above $\beta$ for a given pair $(\alpha,\Delta)$.
The consideration of these questions is motivated by the presence of similar step-by-step procedures in the differential game theory, for example, in the alternating integral method, in pursuit-evasion problems posed with use of counter-strategies, and in the method of guide control.
It is shown that the step-by-step construction of the value $z\in\alpha(\omega)$ can be carried out for any in steps implemented argument $\omega$ if and only if the multiselector $\beta$ is non-empty-valued.
In this case, the desired value $z$ can be selected from $\beta(\omega)$ in step-by-step procedure for any unknown in advance argument $\omega$.
The key point of the work is the procedure for calculation the multiselector $\beta$, for which a constructive and finite-step description is given.
Illustrative examples are considered that include, in particular, problems of a guaranteed result optimization under functional constraints on control and/or disturbance implementations.
}




\maketitle

\section*{Introduction}

Let there be non-empty sets of functions $\naZ$ and $\naOm$ defined on the interval $\naT$.
Let there also be a multifunction (m/f) $S$ defined on $\naOm$ with values in $\naZ$ and a set $N\subset\naZ$.
We will interpret these objects as an abstract control problem under conditions of uncertainty (namely, an abstract game problem of retention):
the uncontrolled factor $\omega$ from the set $\naOm$, acting on the dynamics $S$ of the system, determines the bundle of possible trajectories $S(\omega)\subset\naZ$;
the task of the control side is to select the motion $z\in S(\omega)$ satisfying the phase constraints $N$: $z\in N$.

Without additional informational requirements, the solution of the problem is built explicitly:
we put $\alpha(\omega)\mydef N\cap S(\omega)$, $\omega\in\naOm$;
then the criterion for the solvability of the problem is the non-emptiness of the m/f $\alpha$ values, and the solution is any of its selectors.
So, if we knew the current disturbance $\bar\omega\in\naOm$, the retention problem would be solved by calling any trajectory $z$ of $\alpha(\bar\omega)$.

At the same time, in most control problems, information about the acting disturbance is not available at all.
In the rest cases, the best that the control side can count on by the time $\tau\in\naT$ is the knowledge of the disturbance on the interval $[t_0,\tau+\delta]$ with a small $\delta>0$ ($t_0$ is the initial moment of the process).
Such a prediction of disturbance behavior is admissible in certain control problems, as well as in auxiliary control structures, due to properties of the system dynamics (the mapping $S$) and restrictions describing the set of admissible disturbances $\naOm$.

Thus, we suppose that the control side has the m/f $\alpha$ of global responses and the possibility of a small time-ahead disturbance prediction.
Under these conditions, the following step-by-step scheme of the desired trajectory constructing naturally arises.

Let $\naT\mydef[t_0,\vartheta]$ and $\Delta$, $\Delta\subset\naT$, be a finite set of instants that splits \naT\ into a finite set of half-intervals (control steps) of length shorter then the above value $\delta$:
$\Delta\mydef\{t_0=\tau_0<\tau_1<\ldots<\tau_{n_\Delta}=\vartheta \}$.
We call $\Delta$ \emph{the partition} of $\naT$ and the pair $(\alpha,\Delta)$ will be called \emph{step-by-step conditions}.
The procedure of step-by-step construction of the trajectory that meets the conditions $(\alpha,\Delta)$ and the unknown disturbance $\bar\omega\in\naOm$ works as follows:

--- at the instant $\tau_0$, by the found disturbance $\omega_1\in\naOm$ such that $\res{\omega_1}{[\tau_0,\tau_1)}=\res{\bar\omega}{[\tau_0,\tau_1)}$,
the control side chooses the trajectory $h_1\in\naZ$ corresponding  to the disturbance $\omega_1$, i.e. such that $h_1\in\alpha(\omega_1)$; here, $\res{f}{C}$ is the restriction of a function $f$ to the set $C$;

--- at any instant $\tau_i\in\Delta$, $i\in\nint1{(n_\Delta-1)}$, the control side finds a disturbance $\omega_{i+1}\in \naOm$ reconstructing the unknown $\bar\omega$ up to the moment $\tau_{i+1}$:
\beq\label{omi-rec}
\res{\omega_{i+1}}{[\tau_0,\tau_{i+1})}=\res{\bar\omega}{[\tau_0,\tau_{i+1})};
\eeq
and, so, $\omega_{i+1}$ coincides with $\omega_{i}$ up to the moment $\tau_{i}$
\beq\label{omi-rec-2}
\res{\omega_{i+1}}{[\tau_0,\tau_{i})}=\res{\omega_i}{[\tau_0,\tau_{i})}.
\eeq
At the previous steps of partition $\Delta$, we have already supplied the desired trajectory $h_i$ for the disturbance $\omega_i$, i.e. $h_i\in\alpha(\omega_i)$;
taking this into account, the control side looks for the trajectory $h_{i+1}\in\alpha(\omega_{i+1})$ corresponding to the disturbance $\omega_{i+1}$, that also coincides with the choice at the previous steps:
$$
\res{h_{i+1}}{[\tau_0,\tau_i)}=\res{h_{i}}{[\tau_0,\tau_i)}.
$$
The procedure is repeated for all moments of the partition $\Delta$ except the last one, $\tau_{n_\Delta}$.
As a result, we get the desired trajectory $h_{n_\Delta}$ that corresponds to the unknown in advance disturbance $\omega_{n_\Delta}=\bar\omega$:
$$
h_{n_\Delta}\in\alpha(\omega_{n_\Delta})=\alpha(\bar\omega).
$$
The possibility of realization the above step-by-step procedure in response to any admissible disturbance $\omega\in\naOm$ will be referred to as \emph{the feasibility of the conditions $(\alpha, \Delta)$}.

The work deals with two questions.
The first one is the connection between the feasibility of conditions $(\alpha, \Delta)$ and the existence of some special multi-selector of the m/f $\alpha$.
The second question is the construction of this multiselector for given conditions $(\alpha, \Delta)$.

The motivation for the consideration was the above step-by-step scheme and similar ones, which starting from convergence problem \cite[sec.\,III]{Fleming:1961} in theory of differential games arise, for example, in the method of alternating integral \cite{Pontr1967et}, in pursuit --- evasion problems using counter-strategies  \cite{BlaPet-DGA2019,Chernov2014ARC,PetZen1996et}, or in con\-trol\-ling with a guide under functional constraints on a disturbance (see \cite{CheKhlo2005DU,GomSer-UDSU2021} and references).
Besides game-theoretical problems, the above informational conditions for control side can be found in the field of robotics: suppose a robot-manipulator extracts from the container and submits for the further processing some parts poured into it.
In this case the disturbance/uncertainty is the arrangement of the parts in the container.
It changes (unpredictably as a rule) after the extraction of a part and remains practically constant/unchanged during inactivity time.
So, the disturbance can be effectively predicted.

It is clear that the non-anticipatory and non-emptiness of values of $\alpha$ or of its multiselector implies feasibility of $(\alpha, \Delta)$ for any $\Delta$.
Thus, the conditions for the existence of non-anticipative multiselectors (see, for example, \cite{Chentsov2001DEI, Chentsov2001DEII,Ser_UDSU2017}) are sufficient conditions for the feasibility property.
Note that the existence of a non-anticipative selector \cite{Chentsov1978dep-e,CARDALIAGUET-PLASKACZ-SIAM2000,CheSer_TRIMM2019et} is also closely related to the existence of a non-anticipative multiselector.

In this paper, we show that the feasibility of the conditions $(\alpha, \Delta)$ is equivalent to the existence of a \emph{partially non-anticipative} and non-empty-valued multiselector of the  m/f $\alpha$:
the above step-by-step process can be implemented by means of this  partially non-anticipative multiselector for any disturbance.
Here, the property of partial non-anticipatory is understood as the classical non-anticipatory property that is satisfied at moments from $\Delta$ only.
This property is certainly weaker than the classical one when it should be satisfied for all moments from $\naT$.
Moreover, even the feasibility of the conditions $(\alpha, \Delta)$ for any $\Delta$ does not in general ensure the existence of a non-anticipative non-empty-valued multiselector of $\alpha$.

So, to implement the above step-by-step procedure under conditions $(\alpha, \Delta)$, we need this partially non-anticipative and non-empty-valued multiselector of $\alpha$.
Here arises the second question: how to build such a partially non-anticipative multiselector.

With the aim, for any instant $\tau\in\naT$ we introduce a "projection" operator acting on set of all m/f with values in the set of all $\tau$-non-anticipative m/f.
Then, by means of a  superposition of such "projections" (corresponding to all $\tau\in\Delta$), we get the required partially non-anticipative multiselector.
The procedure is completed in $n_\Delta$ steps.

It seems the idea of constructing a non-anticipative multiselector of a m/f by an iterative method in order to obtain a direct solution of the dynamic optimization problem under conditions of an uncertainty appeared in \cite{Chentsov1997DAN} (see also \cite{Chentsov2001DEI, Chentsov2001DEII}).
The obstacle for applications of such an iterative procedure is that it requires infinite number of iterations in the general case.
For some classes of control problems (see, for example, \cite[ch.\,5]{SubChe81e}), conditions are given that ensure the finiteness of iterations required to construct the function of the optimal guaranteed result, a resolving set of initial positions, or a resolving non-anticipative strategy.

The work is close to the problems considered in \cite{GomSer-UDSU2021}; 
the constructions used are similar to those from \cite{Ser_UDSU2017} and the results supplements the results announced in \cite{Ser_CTMM2022e}.
In simple cases, the connection of partially non-anticipative multiselectors with ordinary non-anticipative multiselectors was noted in the course of the presentation, but is not covered in detail.
The outline of the article is as follows:
Sec.\,\ref{sec-def} contains basic notation and terms;
in Sec.\,\ref{Sbs-cont}, a more detailed and meaningful description of the step-by-step procedure for constructing an optimal trajectory, its formalization, and the definition of the feasibility of the $(\alpha, \Delta)$ conditions are given (see \fref{def-sbs-p1}, \fref{def-sbs-p2});
in Sec.\,\ref{Ana}, the notion of a partially non-anticipative m/f  (multiselector) is defined and the feasibility criterion \fref{aHfiz-eq-aHdOm} is formulated in its terms;
in Sec.\,\ref{Ana-oper}, the above-mentioned "projection" operator on the set of all $\tau$-non-anticipative m/f is defined (see \fref{def-Ana});
then, the definition of such "projection" operator on the set of m/f that are non-anticipative for several such "moments" (see \fref{Ana-max-H}) is given and the finite step procedure for its constructing is provided (see Theorem \ref{fix-to-couple-lem});
in Sec.\,\ref{EXAM}, illustrative examples are given:
in Example \ref{ex2}, it is shown that the "projection" operators are not commutative;
Example \ref{ex3} provides the case when the feasibility of conditions $(\alpha, \Delta)$ for arbitrary $\Delta$ do not allow, nevertheless, a resolving non-anticipative strategy (a non-anticipative non-empty-valued multiselector of $\alpha$);
in Example \ref{ex4}, an ordinal resolving non-anticipative strategy is constructed using the proposed technic.
in Sec.\,\ref{conc}, open formal questions and prospects for using this approach in applications are briefly discussed.

\section{Definitions and notation}
\label{sec-def}

In the following, set-theoretic symbolism is used (quantifiers, propositional connectives, etc.);
hereinafter \myemp\ --- empty set, $\mydef$~--- equality by definition, $\mydefeq$~--- equivalency by definition;
a family is a set all of whose elements are sets.

Let $\icP X$ and $\icPp X$ denote respectively the families of all (Boolean of $X$) and of all non-empty subsets of an arbitrary set $X$.
If $A$ and $B$ are non-empty sets, then we denote by $B^A$ the set of all mappings from $A$ to $B$.
If $g\in\icP B^A$, then by $\nado{g}$ we denote the region where the m/f $g$ takes nonempty values: $\nado{g}\mydef\{a\in A\mid g(a)\neq\myemp\}$.
If $f\in B^A$ and $C\in\icPp A$, then $\res{f}{C}$, where $\res{f}{C}\in B^C$, is the restriction of $f$ to the set $C$: $\res{f}{C}(x)\mydef f(x)$ $\myll x\in C$;
in the case when $F\in\icPp{B^A}$, we set $\sres{F}{C}\mydef\{\res{f}{C}\colon f\in F\}$.

We call a pair $(X,\myles)$ \emph{partially ordered set} (poset) if $X$ is a non-empty set and ${\myles}\in\icP{X\times X}$ is a non-strict partial order relation on $X$.
In particular, $(\icP X,\subset)$ is a poset on Boolean of $X$ with the inclusion relation ${\subset}$.
For any poset $(X,\myles)$ and arbitrary subset $S\in\icP X$ we call $S$ a \emph{chain} if any elements of $S$ are comparable: $(x\myles y)\myor(y\myles x)$, $\myll x,y\in S$.
If $(X,\myles)$ is a poset and $M\subset X$, denote $\top_M$,  $\top_M\in M$, the greatest element of $M$, if it exists: $x\myles\top_M$ for all $x\in M$.

We fix non-empty sets \naX, \naY\ and \naT, as well as non-empty sets $\naOm\in\icPp{\naY^\naT}$, $\naZ\in\icPp{\naX^\naT}$ and the family $\naTc\in\icPp{\icPp\naT}$;
in other words, $\naTc$ is a non-empty family of non-empty subsets of $\naT$.
In order to reveal the connection with dynamic optimization problems, we will assume that $\naT$ is  a "time interval": $\naT\mydef[t_0,\vartheta]$, $t_0,\vartheta\in\RA$ ; and $\naTc$ is a chain of the form $\naTc\mydef\{[t_0,\tau]:\tau\in\naT\}$.
Then $\naOm$ is a set of admissible disturbances, $\naZ$ is a set of possible system trajectories.

Let  $\myLe$ be a partial order on the set mappings $\icP\naZ^\naOm$: $\forall\alpha,\beta\in\icP\naZ^\naOm$
\beq\label{def-Le}
(\alpha\myLe\beta)\mydefeq(\alpha(\omega)\subset\beta(\omega)\ \myll\omega\in\naOm).
\eeq
For $\omega\in\naOm$, $z\in\naZ$ and $A\in\icPp\naT$ we denote
$$
\naOmo\omega A\mydef\{\eta\in\naOm\mid \res{\eta}{A}=\res{\omega}{A}\},\label{Omo-def}\qquad
\naZo zA\mydef\{g\in\naZ\mid \res{g}{A}=\res{z}{A}\}.\label{Zo-def}
$$
Note that the family of partitions $(\{\naOmo\omega A\mid \omega\in\naOm)_{A\in\naTc}$ forms a chain in the set of partitions of $\naOm$ with an embedding relation (see, for example, \cite[Sec.\,3.1]{Engelking1986e}), namely, for any $A,A'\in\naTc$, we have
\beq\label{Omo-isot}
(A\subset A')\myimp((\naOmo\omega{A'}\subset\naOmo\omega A\ \myll\omega\in\naOm).
\eeq
Here by \emph{partition} of a set $M$ we call a family of its nonempty mutually disjont subsets that cover the set.

In terms of the family $\naTc$ we define the basic notion of \emph{non-anticipatory}.
Denote by $\naN$ the set of all \emph{non-anticipative} m/f from ${\icP\naZ}^\naOm$:
$$
\naN\mydef\left\{{\bf z}\in\icP\naZ^\naOm\mid \sres{{\bf z}(\omega)}{A}=\sres{{\bf z}(\omega^ \prime)}{A}\ \myll A\in\naTc\ \myll\omega\in\naOm\ \myll\omega^\prime\in\naOmo\omega A\right\}.
$$
By $\naNO{}$ we denote the subset of all nonempty-valued m/f:
$$
\naNO{}\mydef\{{\bf z}\in\naN\mid\nado{{\bf z}}=\naOm\}.
$$
The most important are multiselectors of a given m/f: for $\alpha\in\icP\naZ^\naOm$ we denote
\beq\label{nano-def}
\naNo\alpha\mydef\{{\bf z}\in\naN\mid{\bf z}\myLe\alpha\};
\eeq
then we define the subset of nonempty-valued multiselectors of m/f $\alpha$:
\beq\label{nanO-def}
\naNO\alpha\mydef\{{\bf z}\in\naNO{}\mid{\bf z}\myLe\alpha\}.
\eeq

Let us introduce (see \cite[\S\,3, item V]{Kura1966e}) pointwise defined "union" and "intersection" of m/f from $\icP\naZ^\naOm$.
Let $\mathbb{M}$ be a subset of $\icP\naZ^\naOm$.
We assume that the pointwise union $\bigvee_{{\bf z}\in\mathbb{M}}{{\bf z}}$ and the pointwise intersection $\bigwedge_{{\bf z}\in\mathbb{M}}{ {\bf z}}$, for each $\omega\in\naOm$ are determined by:
\beq\label{def-mfUI}
\bigvee_{{\bf z}\in\mathbb{M}}{\bf z}(\omega)\mydef\bigcup_{{\bf z}\in\mathbb{M}}{{\bf z}} (\omega),\qquad
\bigwedge_{{\bf z}\in\mathbb{M}}{\bf z}(\omega)\mydef\bigcap_{{\bf z}\in\mathbb{M}}{{\bf z}} (\omega).
\eeq

Note that in the form $\bigvee_{{\bf z}\in\mathbb{M}}{\bf z}$ and $\bigwedge_{{\bf z}\in\mathbb{M}}{\bf z }$ we have, respectively, supremum and  infimum of the set $\mathbb{M}$ in the poset $(\icP\naZ^\naOm,\myLe)$.
Unlike pointwise intersection, the result of pointwise union inherits non-anticipatory property of operands (see lemma \ref{lem-AnaUI-Ana} below).

\section{Feasibility of step-by-step procedure}
\label{Sbs-cont}

Let us define the property of feasibility in formal terms.
Let $\Delta$, $\Delta\subset\naT$, be a finite set of instants: $\Delta\mydef\{t_0=\tau_0<\tau_1<\ldots<\tau_{n_\Delta}=\vartheta \}$.
We will also refer to $\Delta$ as the partition of the time interval \naT.
Denote by ${\cal H}_\Delta$ the subset of $\naTc$ of the form ${\cal H}_\Delta\mydef\{H_i\mydef[\tau_0,\tau_i]: i\in\nint1{n_\Delta}\}$.

For partition $\Delta$ and a m/f $\alpha\in\icP\naZ^\naOm$ denote by $\naOm_\Delta$ ($\naOm_\Delta\subset\naOm^{n_\Delta}$) and $\naZ_\Delta$ ($\naZ_\Delta\subset{\icP\naZ}^{n_ \Delta}$) the sets defined as follows (in expressions like $Q^n$, where $n\in\NA$, the number $n$ is not considered as a set):
\beq\label{OD-def-2}
\naOm_\Delta\mydef\{(\omega_i)_{i\in\nint1{n_\Delta}}\in\naOm^{n_\Delta}\mid\res{\omega_i}{H_{i}}=\res{\omega_{i+1}}{H_{i}},\ i\in\nint1{n_\Delta-1}\},
\eeq
\beq\label{ZD-def-2}
\naZ_\Delta\mydef\{({\beta}_i)_{i\in\nint1{n_\Delta}}\in{\icPp\naZ}^{n_\Delta}\mid\sres{{\beta}_i}{H_{i}}=\sres{{\beta}_{i+1}}{H_{i}},\ i\in\nint1{n_\Delta -1}\}.
\eeq
Conditions $(\alpha, \Delta)$ for the step-by-step procedure we call feasible if there exists a tuple of multiselectors $(\phi_i)_{i\in\nint1{n_\Delta}}$ of  m/f $\alpha$,
\beq\label{def-sbs-p1}
\phi_i\in\icP\naZ^\naOm,\quad\phi_i\myLe\alpha,\qquad i\in\nint1{n_\Delta},
\eeq
such that the following inclusions are fulfilled:
\beq\label{def-sbs-p2}
(\phi_i(\omega_i))_{i\in\nint1{n_\Delta}}\in\naZ_\Delta\qquad\forall (\omega_i)_{i\in\nint1{n_\Delta}}\in\naOm_\Delta.
\eeq

This definition simply retells the above informal description of the step-by-step procedure in mathematical terms.
Namely, let given a tuple $(\phi_i)_{i\in\nint1{n_\Delta}}$ satisfying \fref{def-sbs-p1}, \fref{def-sbs-p2}, then the construction of a step-by-step response to an unknown disturbance $\bar\omega$, due to the $(\alpha,\Delta)$ conditions, can be realised like this:

--- we find $\omega_1\in\naOm$ such that $\res{\omega_1}{H_1}=\res{\bar\omega}{H_1} $ and choose $h_1\in\naZ$ that satisfies $h_1\in\phi_1(\omega_1)$.
This can always be done, because by completing formally $\omega_1$ to an arbitrary tuple $(\omega_1,\ldots,\omega_{n_\Delta})$ from $\naOm_\Delta$ (for example, by setting $\omega_{n_\Delta}\mydef\ldots\mydef\omega_2\mydef\omega_1$), due to \fref{ZD-def-2} and \fref{def-sbs-p2} we get the inequality $\phi_1(\omega_1)\neq\varnothing$ which enables us to choose the desired trajectory $h_1$;

--- if there are $(\omega_i)_{i\in\nint1{k+1}}\in\naOm^{k+1}$ and $h_k\in\naZ$ such that
\begin{gather}
\res{\omega_{k+1}}{H_{k+1}}=\res{\bar\omega}{H_{k+1}} ,\label{bomip1}\\
\res{\omega_{i+1}}{H_{i}}=\res{\omega_i}{H_{i}},\qquad i\in\nint1k\label{bomip12}\\
h_i\in\phi_i(\omega_i),\qquad i\in\nint1k\label{hiin}
\end{gather}
we choose $h_{k+1}$ from the condition
$$
h_{k+1}\in\phi_{k+1}(\omega_{k+1})\cap\naZo{h_k}{H_k}.
$$
Such a choice is possible: the tuple $(\omega_1,\ldots,\omega_{k+1})$ can be completed (see \fref{bomip12}, \fref{OD-def-2}) to a tuple $(\omega_1,\ldots,\omega_{n_\Delta})$ from $\naOm_\Delta$ (for example, by setting $\omega_{n_\Delta} \mydef \ldots \mydef \omega_{k+2} \mydef \omega_{k+1}$).
Then, by virtue of \fref{ZD-def-2}, \fref{def-sbs-p2} we get the equality
$$
\res{\phi_{k+1}(\omega_{k+1})}{H_{k}}=\res{\phi_{k}(\omega_{k})}{H_{k}},
$$
from which, taking into account \fref{hiin}, we get $\phi(\omega_{k+1})\cap\naZo{h_k}{H_k}\neq\varnothing$.
Hence, the choice is possible to undertake.
So, by induction, the step-by-step procedure can be continued up to $k+1=n_\Delta$.

It is clear that the non-emptiness of values and non-anticipatory property of $\alpha$ (that is, when $\alpha\in\naNO{}$) or, in general, the existence of a non-empty-valued non-anticipative multiselector $\beta$ of m/f $\alpha$ implies the feasibility of the conditions $(\alpha, \Delta)$ for any partition $\Delta$.
Indeed, then it suffices to put $\phi_i\mydef\beta$, where $\beta\in\naNO\alpha$.
At the same time, it can be seen from the example \ref{ex3} that the feasibility of the conditions $(\alpha, \Delta)$ even for all partitions $\Delta$ do not imply existence of non-anticipative and non-empty-valued multiselector of $\alpha$.
In following sections, we consider the property of partial non-anticipatory of a m/f that is equivalent to the feasibility property.

\section{Partially non-anticipative mappings: basic properties}
\label{Ana}

For an arbitrary $A\in\naTc$, we introduce the notion of \emph{$A$-non-anticipative} m/f:
a mapping $\mathbf{z}\in{\icP\naZ}^\naOm$ is called $A$-non-anticipative if the implication
$$
\left(\res{\omega_1}{A}=\res{\omega_2}{A}\right) \myimp \left(\sres{\mathbf{z}(\omega_1)}{A}=\sres{\mathbf{z}(\omega_2)}{A}\right)
$$
is true for all $\omega_1,\omega_2\in\naOm$.
The subset of all $A$-non-anticipative m/f we denote by $\AnaN{\{A\}}{}$.
The subset of all non-empty-valued and $A$-non-anticipative m/f we denote by $\AnaNO{\{A\}}{}$:
$\AnaNO{\{A\}}{ }\mydef\{{\bf z}\in\AnaN{\{A\}}{}\mid \nado{\bf z}=\naOm\}$.

Let ${\cal H}\in\icPp\naTc$.
We call a m/f from ${\icP\naZ}^\naOm$ \emph{${\cal H}$-non-anticipative} if the m/f is $A$-non-anticipative for all $A\in{\cal H}$.
The set all ${\cal H}$-non-anticipative (non-empty-valued and ${\cal H}$-non-anticipative) m/f we denote by $\AnaN{\cal H}{ }$ ($\AnaNO{\cal H}{}$):
$$
\AnaN{\cal H}{}=\bigcap_{A\in{\cal H}}\AnaN{\{A\}}{}\qquad\bigl( \AnaNO{\cal H}{}=\bigcap_{A\in{\cal H}}\AnaNO{\{A\}}{}\bigr).
$$
For any $\alpha\in\icP\naZ^\naOm$ and ${\cal H}\in\icPp\naTc$ designations $\AnaN{\cal H}\alpha$ and $\AnaNO{\cal H}\alpha$  are defined similarly to those in \fref{nano-def} and \fref{nanO-def}, respectively.

Note the implication following from these definitions (isotonicity):
for arbitrary $\alpha\in\icP\naZ^\naOm$, ${\cal H}$, ${\cal H'}\in\icP\naTc$
\beq\label{A-nan-isot-H}
\left({\cal H}\subset{\cal H'}\right)\myimp\left(\AnaNO{\cal H'}{\alpha}\subset\AnaNO{\cal H}{\alpha}\right).
\eeq

The following is definitions and some properties of point-wise operations on non-anticipative m/f.
\begin{lem}\label{lem-AnaUI-Ana}
Let $\alpha\in\icP\naZ^\naOm$, ${\cal H}\subset\naTc$ and $\mathbb{M}\in\icPp{\AnaN{\cal H}\alpha}$.
Then
\beq\label{AnaUI-Ana}
\bigvee_{{\bf z}\in\mathbb{M}}{\bf z}\in\AnaN{\cal H}\alpha,
\eeq
\beq\label{zom-nem-ANo-nem}
\Bigl(\myll\omega\in\naOm\ \myxst{\bf z}\in\mathbb{M}: {\bf z}(\omega)\neq\myemp\Bigr)\myimp\left(\bigvee_ {{\bf z}\in\mathbb{M}}{\bf z}\in\AnaNO{\cal H}\alpha\right),
\eeq
\beq\label{def-Ana-al}
\top_{\AnaN{\cal H}\alpha}=\bigvee_{{\bf z}\in\AnaN{\cal H}\alpha}{\bf z}.
\eeq
\end{lem}

\proof
1.
It is clear that
\beq\label{AnaUI-Ana1}
\left(\bigvee_{{\bf z}\in\mathbb{M}}{\bf z}\myLe\alpha\right)\myand\left(\bigwedge_{{\bf z}\in\mathbb{M }}{\bf z}\myLe\alpha\right).
\eeq
We show that the mapping $\bigvee_{{\bf z}\in\mathbb{M}}{\bf z}$ inherits the property of ${\cal H}$-non-anticipatory of elements from $\mathbb{M}$: let $ \omega,\omega'\in\naOm$, $A\in{\cal H}$, $\res\omega A=\res{\omega'}A$ and $\gamma\in\sres{\bigvee_ {{\bf z}\in\mathbb{M}}{\bf z}(\omega)}A$.
Then there are $\bar{\bf z}\in\mathbb{M}$ and $h\in\bar{\bf z}(\omega)$ such that $\gamma=\res hA$.
Therefore, taking into account the relation $\bar{\bf z}\in\AnaN{\cal H}\alpha$, we have the inclusion $\gamma\in\sres{\bar\bf z(\omega')}A$.
Hence, there exists $h'\in\bar{\bf z}(\omega')$ such that $\gamma=\res{h'}A$, whence we obtain the inclusion $h'\in\bigvee_{{\bf z}\in\mathbb{M}}{\bf z}(\omega')$ and, as a consequence, the inclusion $\gamma\in\sres{\bigvee_{{\bf z}\in\mathbb{M }}{\bf z}(\omega')}A$.
Since $\gamma$ was chosen arbitrarily, we have the relation
$$
\sres{\bigvee_{{\bf z}\in\mathbb{M}}{\bf z}(\omega)}A\subset\sres{\bigvee_{{\bf z}\in\mathbb{M} }{\bf z}(\omega')}A.
$$
Hence, due to symmetry of considerations and an arbitrary choice of $\omega$, $\omega'$, the implication follows
$$
(\res\omega A=\res{\omega'} A)\myimp\sres{\bigvee_{{\bf z}\in\mathbb{M}}{\bf z}(\omega)}A=\sres{\bigvee_{{\bf z}\in\mathbb{M}}{\bf z}(\omega')}A\qquad\forall\omega,\omega'\in\naOm.
$$
Due to the implication, taking into account the arbitrary choice of $A$, we have
\beq\label{AnaUI-Ana2}
\bigvee_{{\bf z}\in\mathbb{M}}{\bf z}\in\AnaN{\cal H}{}.
\eeq
Relations \fref{AnaUI-Ana1}, \fref{AnaUI-Ana2} together give the desired inclusion \fref{AnaUI-Ana}.

2.
Premise in \fref{zom-nem-ANo-nem} implies the equality $\nado{\bigvee_{{\bf z}\in\mathbb{M}}{\bf z}}=\naOm$, which in combination with \fref{AnaUI-Ana} gives the desired assertion.

3.
By construction, the expression on the right side \fref{def-Ana-al} $\myLe$-majorizes each element of the set $\AnaN{\cal H}\alpha$, while \fref{AnaUI-Ana} for $\mathbb{ M}\mydef\AnaN{\cal H}\alpha$ implies the inclusion $\bigvee_{{\bf z}\in\AnaN{\cal H}\alpha}{\bf z}\in\AnaN{\cal H }\alpha$.

The lemma is proven.

\bigskip

\begin{theorem}
\label{lem-aHfiz-eq-aHdOm}
For a m/f $\alpha\in\icP\naZ^\naOm$ and a partition $\Delta$ of interval $\naT$, the conditions $(\alpha,\Delta)$ are feasible if and only if m/f $\alpha$ has an non-empty-valued and ${\cal H}_\Delta$-non-anticipative multiselector:
\beq\label{aHfiz-eq-aHdOm}
\left((\alpha,\Delta)-\text{feasible}\right)\myeqv\left(\AnaNO{{\cal H}_\Delta}\alpha\neq\varnothing\right).
\eeq
\end{theorem}
\proof
Let $\Delta\mydef\{t_0=\tau_0<\tau_1<\ldots<\tau_{n_ \Delta}=\vartheta\}$.
Remind that ${\cal H}_\Delta\mydef\{H_i\mydef[\tau_0,\tau_i]: i\in\nint1{n_\Delta}\}$.
The definitions \fref{OD-def-2}, \fref{ZD-def-2} imply the equalities
\beq\label{OD-def}
\naOm_\Delta=\{(\omega_i)_{i\in\nint1{n_\Delta}}\in\naOm^{n_\Delta}\mid\res{\omega_i}{H_{ \overline{ij}}}=\res{\omega_j}{H_{\overline{ij}}},\ i,j\in\nint1{n_\Delta}\},
\eeq
\beq\label{ZD-def}
\naZ_\Delta=\{({\bf z}_i)_{i\in\nint1{n_\Delta}}\in{\icPp\naZ}^{n_\Delta}\mid\sres{{\bf z}_i}{H_{\overline{ij}}}=\sres{{\bf z}_j}{H_{\overline{ij}}},\ i ,j\in\nint1{n_\Delta}\},
\eeq
where $\overline{ij}\mydef\min\{i,j\}$ for all $i,j\in\nint1{n_\Delta}$.

1.
Let $\AnaNO{{\cal H}_\Delta}\alpha\neq\varnothing$ and $\beta\in\AnaNO{{\cal H}_\Delta}\alpha$.
Consider a tuple $(\phi_i)_{i\in\nint1{n_\Delta}}$ of the form $\phi_i\mydef\beta$, $i\in\nint1{n_\Delta}$.
By definition, it satisfies the equalities $\nado{\phi_i}=\naOm$, $i\in\nint1{n_\Delta}$, and the conditions $\phi_i\myLe\beta\myLe\alpha$, $i\in\nint1{n_ \Delta}$.
This implies the fulfillment of the condition \fref{def-sbs-p1}.

Let us check \fref{def-sbs-p2}.
Suppose a tuple $(\omega_i)_{i\in\nint1{n_\Delta}}$ is such that $(\omega_i)_{i\in\nint1{n_\Delta}}\in\naOm_\Delta$.
Then, taking into account the ${\cal H}_\Delta$-non-anticipatory property of $\beta$, from the equalities $\res{\omega_k}{H_{\overline{km}}} =\res{\omega_m}{H_{\overline{km}}}$ we get the equalities
$$
\sres{\phi_k(\omega_k)}{H_{\overline{km}}} =\sres{\beta(\omega_k)}{H_{\overline{km}}} =\sres{\beta(\omega_m)}{H_{\overline{km}}} =\sres{\phi_m(\omega_m)}{H_{\overline{km}}}\quad k,m\in\nint1{n_\Delta}.
$$
Since $(\omega_i)_{i\in\nint1{n_\Delta}}$, $k$ and $m$ were chosen arbitrarily, for the tuple $(\phi_i)_{i\in\nint1{n_\Delta }}$ condition \fref{def-sbs-p2} is met.
We have shown that the left side of \fref{aHfiz-eq-aHdOm} follows from the right side.

2.
Let us show that the right side of \fref{aHfiz-eq-aHdOm} follows from the left side.
Assume that the conditions $(\alpha,\Delta)$ are feasible.
Then (see sec.\,\ref{Sbs-cont}) there exists a tuple $(\phi_i)_{i\in\nint1{n_\Delta}}$ of the form \fref{def-sbs-p1} such that for any tuple $(\omega_i)_{i\in\nint1{n_\Delta}}\in\naOm_\Delta$ (see \fref{OD-def},\fref{ZD-def}):
\beq\label{fia-in-ZD}
(\phi_i(\omega_i))_{i\in\nint1{n_\Delta}}\in\naZ_\Delta.
\eeq
Let m/f $\phi_\alpha$ be defined by $\phi_\alpha\mydef\phi_{n_\Delta}$.
Then we have the comparison $\phi_\alpha\myLe\alpha$ and  equality $\nado{\phi_\alpha}=\naOm$
(indeed, from \fref{ZD-def}, \fref{fia-in-ZD} follow the inclusions $\phi_\alpha(\omega)\in\icPp\naZ$ for all $\omega\in\naOm$).
To show that $\phi_\alpha$ is ${\cal H}_\Delta$-non-anticipative, suppose $\omega,\omega'\in\naOm$ and $m\in\nint1{n_\Delta}$ are such that
\beq\label{om-ta-hat}
\res\omega{H_{m}}=\res{\omega'}{H_{m}}.
\eeq
We put
$$
\omega_i\mydef\omega,\quad i\in\nint1{n_\Delta},\qquad
\omega'_i\mydef
\begin{cases}
\omega,& i\in\nint1{m},\\
\omega',& i\in\nint{(m+1)}{n_\Delta}.
\end{cases}
$$

For the tuple $(\omega_i)_{i\in\nint1{n_\Delta}}$, we obviously have the inclusion $(\omega_i)_{i\in\nint1{n_\Delta}}\in\naOm_\Delta $.

Show the inclusion $(\omega'_i)_{i\in\nint1{n_\Delta}}\in\naOm_\Delta$.
For the tuple $(\omega'_i)_{i\in\nint1{n_\Delta}}$ and any $i,j\in\nint1{n_\Delta}$ we have:

if $m< i\le j$, then
$$
\res{\omega'_i}{H_{\overline{ij}}}=\res{\omega'}{H_{\overline{ij}}}=\res{\omega '_j}{H_{\overline{ij}}};
$$

if $i\le m<j$, then (see \fref{om-ta-hat})
$$
\res{\omega'_i}{H_{\overline{ij}}}=\res{\omega}{H_{\overline{ij}}}=\res{\omega}{H_{i}}=\res{\omega' }{H_{i}}=\res{\omega'_j }{H_{i}}=\res{\omega'_j}{H_{\overline{ij}}};
$$

if $i\le j\le m$, then
$$
\res{\omega'_i}{H_{\overline{ij}}}=\res{\omega}{H_{\overline{ij}}}=\res{\omega' _j}{H_{\overline{ij}}}.
$$
Thus (see \fref{OD-def}), the inclusion $(\omega'_i)_{i\in\nint1{n_\Delta}}\in\naOm_\Delta$ takes place.

Therefore, due to \fref{fia-in-ZD}, the inclusions
$$
(\phi_i(\omega_i))_{i\in\nint1{n_\Delta}}\in\naZ_\Delta,\qquad (\phi_i(\omega'_i))_{i\in\nint1{n_\Delta}}\in\naZ_\Delta,
$$
considered (see \fref{ZD-def}) under $i=m$ and $j=n_\Delta$, imply the equalities
\beq\label{fmom-eq-fnon}
\sres{\phi_{m}(\omega_{m})}{H_{\overline{m n_\Delta}}}=\sres{\phi_{n_\Delta}(\omega_{ n_\Delta})}{H_{\overline{m n_\Delta}}},\qquad
\sres{\phi_{m}(\omega'_{m})}{H_{\overline{m n_\Delta}}}=\sres{\phi_{n_\Delta}(\omega'_{n_\Delta})}{H_{\overline{m n_\Delta}}}.
\eeq
Moreover, from the definitions of $\omega'_{m}$ and $\omega_{m}$ we have $\omega'_{m}=\omega_{m}=\omega$ and hence the equalities
\beq\label{pmom-eq-pmopm}
\sres{\phi_{m}(\omega_{m})}{H_{\overline{m n_\Delta}}}=\sres{\phi_{m}(\omega)}{H_{m}}=\sres{\phi_{m}(\omega'_{m})}{H_{\overline{m n_\Delta}}}.
\eeq
From the given equalities and the definition of $\phi_\alpha$ we obtain (the 2nd and the 4th equalities follow from \fref{ZD-def} and \fref{fmom-eq-fnon}, the 3rd --- from \fref{pmom-eq-pmopm}):
\begin{multline*}
\sres{\phi_\alpha(\omega)}{[\tau_0,\tau_{m}]}=\sres{\phi_{n_\Delta}(\omega_{n_\Delta})}{[\tau_0, \tau_{\overline{m n_\Delta}}]}=\sres{\phi_{m}(\omega_{m})}{[\tau_0,\tau_{\overline{m n_\Delta}}]} \\
=\sres{\phi_{m}(\omega'_{m})}{[\tau_0,\tau_{\overline{m n_\Delta}}]}=\sres{\phi_{n_\Delta}( \omega'_{n_\Delta})}{[\tau_0,\tau_{\overline{m n_\Delta}}]}=\sres{\phi_\alpha(\omega')}{[\tau_0,\tau_{m}]}.
\end{multline*}
Since $m$, $\omega$ and $\omega'$ were chosen arbitrarily, from the last equalities we obtain the property of ${\cal H}_\Delta$-non-anticipatory of m/f $\phi_\alpha$.
Taking into account the indicated properties of $\phi_\alpha$, we have the inclusion $\phi_\alpha\in\AnaNO{{\cal H}_\Delta}{\alpha}$, i.e., the right side of \fref{aHfiz-eq-aHdOm} is fulfilled.

The proof is complete.

\section{A construction of the partially non-anticipative multiselector}
\label{Ana-oper}

In this section, we give a description of partially non-anticipative multiselectors in terms of explicitly defined operators that are non-anticipative at one point of $\naT$.
In general, such description is certainly non-constructive.
Meanwhile, when applied to step-by-step procedures, due to the finiteness of operations, the description allows to construct and analyze corresponding partially non-an\-ti\-ci\-pat\-ive multiselectors (see examples).

For an arbitrary $A\in\naTc$, denote by $\Ana{\cdot}{A}$ the operator that transforms the set $\icP\naZ^\naOm$ and is given by:
\beq\label{def-Ana}
\Ana{\alpha}{A}(\omega)\mydef\bigl\{h\in\alpha(\omega)\ \bigr|\ \res hA\in\bigcap_{\omega'\in\naOmo\omega A}\sres{\alpha(\omega')}A\bigr\}\qquad \forall\alpha\in\icP\naZ^\naOm,\ \forall\omega\in\naOm.
\eeq

It immediately follows that  $\Ana{\cdot}{A}$ is non-expansive and isotonic as an operator in the poset $(\icP\naZ^\naOm, \myLe)$:
for arbitrary $A\in\naTc$, $\alpha,\beta\in\icP\naZ^\naOm$
\beq\label{Ana-cont}
\Ana{\alpha}{A}\myLe\alpha,
\eeq
\beq\label{Ana-isot}
(\alpha\myLe\beta)\Rightarrow\left(\Ana{\alpha}{A}\myLe\Ana{\beta}{A}\right).
\eeq
We also note (see subsection \ref{ex1}) that in the general case the mapping $\naTc\ni A\mapsto\Ana{\alpha}{A}\in\icP\naZ^\naOm$ is not isotonic as a mapping from poset $(\naTc,\subset)$ to poset $(\icP\naZ^\naOm,\myLe)$.

\begin{lem}\label{lem-Ana-max}
For $A\in\naTc$ and $\alpha\in\icP\naZ^\naOm$

{\rm (i)} values of the operator $\Ana{\cdot}{A}$ consist the set of all $A$-non-anticipative m/f, as well as the set of its fixed points:
\beq\label{Ana-fix}
\Ana{\icP\naZ^\naOm}{A}=\AnaN{\{A\}}{}=\fix{\Ana{\cdot}{A}};
\eeq

{\rm (ii)} $\Ana{\alpha}{A}$ is the $\myLe$-greatest $A$-non-anticipative multiselector of $\alpha$:
\beq\label{Ana-max}
\Ana{\alpha}{A}=\top_{\AnaN{\{A\}}\alpha};
\eeq

{\rm (iii)} the operator \Ana{\cdot}A is idempotent:
for any $\alpha\in\icP\naZ^\naOm$ the following equality takes place
\beq\label{Ana-idemp}
\Ana{\Ana{\alpha}A}A=\Ana{\alpha}A.
\eeq
\end{lem}

\proof
1.
Let $\beta\mydef\Ana{\alpha}{A}$.
We show that $\beta\in\AnaN{\{A\}}\alpha$.
Due to \fref{Ana-cont}, we have $\beta\myLe\alpha$.
It remains to verify the $A$-non-anticipatory property of $\beta$.

Let $\omega,\omega'\in\naOm$ be such that $\res\omega A=\res{\omega'}A$ and $\xi\in\sres{\beta(\omega)}A $.
Then (see \fref{def-Ana})
\beq\label{xi-in-bigcap}
\xi\in\bigcap_{\bar\omega\in\naOmo\omega A}\sres{\alpha(\bar\omega)}A.
\eeq
By the choice of $\omega'$, we have $\omega'\in\naOmo\omega A$ and hence (see \fref{xi-in-bigcap}) $\xi\in\sres{\alpha(\omega')}A$.
Then, there exists $h'\in\alpha(\omega')$ such that $\res{h'}A=\xi$.
From the equality $\res\omega A=\res{\omega'}A$ it also follows that $\naOmo{\omega'} A=\naOmo\omega A$ and, therefore,
$$
\bigcap_{\bar\omega\in\naOmo\omega A}\sres{\alpha(\bar\omega)}A=\bigcap_{\bar\omega\in\naOmo{\omega'} A}\sres{ \alpha(\bar\omega)}A.
$$
As a result, $h'$ satisfies relations
$$
h'\in\alpha(\omega'),\qquad\res{h'}A\in\bigcap_{\bar\omega\in\naOmo{\omega'} A}\sres{\alpha(\bar\omega)}A,
$$
i.e. (see \fref{def-Ana}) $h'\in\Ana{\alpha}{A}(\omega')$.
Then $\xi\in\sres{\beta(\omega')}A$.
Since $\xi$ was chosen arbitrarily, we have
$$
\sres{\beta(\omega)}A\subset\sres{\beta(\omega')}A.
$$
From the inclusion, due to an arbitrary choice of $\omega$, $\omega'$ and to symmetry of them in the considerations, we obtain the desired $A$-non-anticipatory of $\beta$.

So, taking into account the definition of $\beta$, we have the inclusions
\beq\label{ana-in-anan}
\Ana{\alpha}{A}\in\AnaN{\{A\}}\alpha\subset\AnaN{\{A\}}{}.
\eeq
Due to the arbitrary choice of $\alpha$, \fref{ana-in-anan} implies an embedding
\beq\label{anaZO-in-NA}
\Ana{\icP\naZ^\naOm}{A}\subset\AnaN{\{A\}}{}.
\eeq
2.
Let us verify that $\Ana{\alpha}{A}$ is the $\myLe$-greatest m/f in $\AnaN{\{A\}}\alpha$.
Let
\beq\label{Ana-max1}
\beta\in\AnaN{\{A\}}\alpha,
\eeq
$\bar\omega\in\naOm$ and $\bar h\in\beta(\bar\omega)$.
Then \fref{Ana-max1} implies that
$$
\res{\bar h}A\in\sres{\beta(\omega)}A\subset\sres{\alpha(\omega)}A\quad\forall\omega\in\naOmo{\bar\omega} A.
$$
So, $\bar h$ satisfies the relations
$$
\bar h\in\beta(\bar\omega)\subset\alpha(\bar\omega),\qquad\res{\bar h}A\in\bigcap_{\omega\in\naOmo{\bar\omega } A}\sres{\beta(\omega)}A\subset\bigcap_{\omega\in\naOmo{\bar\omega } A}\sres{\alpha(\omega)}A,
$$
i.e. (see \fref{def-Ana}) $\bar h\in\Ana{\beta}{A}(\bar\omega)$ and $\bar h\in\Ana{\alpha}{A }(\bar\omega)$.
Since $\bar\omega$ and $\bar h$ were chosen arbitrarily, we have comparisons
\beq\label{bet-in-ana-bet}
\beta\myLe\Ana{\beta}{A},
\eeq
\beq\label{bet-in-ana-al}
\beta\myLe\Ana{\alpha}{A}.
\eeq
Comparing \fref{bet-in-ana-al}, due to arbitrary choice of $\beta$, gives the equality \fref{Ana-max}.

3.
Relation \fref{bet-in-ana-bet} and the non-expansion property of $\Ana{\cdot}{A}$ (see \fref{Ana-cont}) imply the equality $\beta=\Ana{\beta} {A}$.
Hence, due to arbitrary choice of $\beta$, we obtain for any $\alpha\in\icP\naZ^\naOm$ the embedding $\AnaN{\{A\}}{\alpha}\subset\fix{\Ana{\cdot}{A}}$.
For $\alpha$ of the form $\alpha(\omega)\mydef\naZ$, $\omega\in\naOm$, we obviously have $\AnaN{\{A\}}\alpha=\AnaN{\{A \}}{}$.
So the following inclusion is correct
\beq\label{NA-al-in-fix-ana-al}
\AnaN{\{A\}}{}\subset\fix{\Ana{\cdot}{A}}.
\eeq
For any $\beta\in\fix{\Ana{\cdot}{A}}$,  by the definition of a fixed point the equality $\beta=\Ana{\beta}{A}$ is satisfied, i.e., $\beta$ lies in the image of the set $\icP\naZ^\naOm$ under the mapping $\Ana{\cdot}{A}$.
Then the next inclusion is fulfilled:
\beq\label{fix-ana-in-anaZO}
\fix{\Ana{\cdot}{A}}\subset\Ana{\icP\naZ^\naOm}{A}.
\eeq
From relations \fref{anaZO-in-NA}, \fref{NA-al-in-fix-ana-al} and \fref{fix-ana-in-anaZO}, we get the equalities \fref{Ana-fix}.

4.
Equalities \fref{Ana-fix} imply the equality \fref{Ana-idemp}.

The lemma is proven.

\bigskip

From the lemma, we immediately obtain a corollary that allows us to filter out m/f that do not have an non-empty-valued and non-anticipative multiselector.

\begin{cor}\label{lem-Ana-max-est}
For $\alpha\in\icP\naZ^\naOm$ and ${\cal H}\in\icPp\naTc$ the following implications hold:
\beq\label{Ana-max-est}
\left(\beta\in\AnaN{\cal H}\alpha\right)\Rightarrow\left(\beta\myLe\bigwedge_{A\in{\cal H}}\Ana{\alpha}{A}\right),
\eeq
\beq\label{AnaH-na-Om}
\left(\nado{\bigwedge_{A\in{\cal H}}\Ana{\alpha}{A}}\neq\naOm\right)\Rightarrow\left(\AnaNO{\cal H}\alpha= \varnothing\right)\Rightarrow\left(\naNO\alpha=\varnothing\right).
\eeq
\end{cor}

\proof
1.
Let $\beta$ satisfy the premise of \fref{Ana-max-est}.
Then from the relations
$$
\beta\in\AnaN{\{A\}}\alpha, \qquad\forall A\in{\cal H},
$$
and \fref{Ana-max} we have
$$
\beta\myLe\Ana{\alpha}{A}, \qquad\forall A\in{\cal H},
$$
which implies the conclusion of \fref{Ana-max-est}.

2.
From \fref{Ana-max-est} and \fref{A-nan-isot-H} the implications \fref{AnaH-na-Om} follow.

The the corollary is proven.

\bigskip

Since for arbitrary ${\cal H}\in\icPp\naTc$ and $\alpha\in\icP\naZ^\naOm$ there exists the greatest (unique) element $\top_{\AnaN{\cal H} {\alpha}}$ (see the lemma \ref{lem-AnaUI-Ana}), then we introduce an operator $\Ana{\cdot}{{\cal H}}:\icP\naZ^\naOm\mapsto \icP\naZ^\naOm$ of the form
\beq\label{Ana-max-H}
\Ana{\alpha}{{\cal H}}\mydef\top_{\AnaN{\cal H}\alpha},\qquad\alpha\in\icP\naZ^\naOm.
\eeq
For ${\cal H}\subset\naTc$ we denote by ${\bf F}_{\cal H}$ the family of operators defined by ${\bf F}_{\cal H}\mydef\{\Ana{\cdot}{A}\mid A\in{\cal H}\}$.
By $\fix{{\bf F}_{\cal H}}$ we refer to the set of joint fixed points of the family ${\bf F}_{\cal H}$: $\fix{{\bf F}_{\cal H}}\mydef\bigcap_{A\in{\cal H}}\fix{\Ana{\cdot}{A}}$.

\begin{remark}\label{zam-Ana-max_H}
It follows from the definition of $\AnaN{\cal H}{}$ and the equality \fref{Ana-fix} that the joint fixed points of the family ${\bf F}_{\cal H}$ are ${\cal H}$-non-anticipative mappings: $\fix{{\bf F}_{\cal H}}=\AnaN{\cal H}{}$.
Thus, in the case when ${\cal H}$ is a singleton (${\cal H}=\{A\}$, $A\in\naTc$), we get the equality (see \fref{Ana-max}) $\Ana{\cdot}{\{A\}}=\Ana{\cdot}A$.
\end{remark}

\bigskip

We turn to the representation of the operator $\Ana{\cdot}{\cal H}$ in terms of superposition of operators from ${\bf F}_{\cal H}$.
As already was noted, the representation of non-anticipative multiselectors in the form of fixed points of non-expansive isotonic operators in a poset and, as a consequence, in the form of limits of their iterative sequences, was proposed and studied in \cite{Chentsov1997DAN,Chentsov2001DEI,Chentsov2001DEII}.
The iterative process in some cases turns out to be finite (see, for example, \cite[ch.\,5]{SubChe81e}), which makes it possible to obtain efficient solutions to problems.

In this paper despite the fact that the operators from ${\bf F}_{\cal H}$ are in general non-commutative (see subsection \ref{ex2}), a finite-step construction of the operator $\Ana{\cdot}{\cal H}$ is given in theorem \ref{fix-to-couple-lem} for all cases when the set ${\cal H}$ is finite.
To this end, we give some definitions and auxiliary results.

\begin{lem}\label{utv1}
Let $H_1,H_2\in\naTc$, $\alpha\in\icP{\naZ}^\naOm$, $\Upsilon\in\icPp{\naOm}$ and $K\in\icP{\sres{ \naZ}{H_1}}$ are such that
\begin{gather}
H_1\subset H_2,\label{utv-e0}\\
\sres{\alpha(\omega)}{H_2}=\sres{\alpha(\omega')}{H_2},\qquad\forall\omega,\omega'\in\Upsilon.\label{utv-e1}
\end{gather}
Then for $\beta\in\icP{\naZ}^\naOm$ of the form
\beq\label{utv1-beta}
\beta(\omega)\mydef
\begin{cases}
\{h\in\alpha(\omega)\mid\res{h}{H_1}\in K\}, &\omega\in\Upsilon,\\
\alpha(\omega), &\omega\in\naOm\setminus\Upsilon,
\end{cases}
\eeq
the equalities \fref{utv-e1} are also satisfied:
\beq\label{utv-e2}
\sres{\beta(\omega)}{H_2}=\sres{\beta(\omega')}{H_2},\qquad\forall\omega,\omega'\in\Upsilon.
\eeq
\end{lem}

\proof
Let $\omega,\omega'\in\Upsilon$ and $\gamma\in\sres{\beta(\omega)}{H_2}$.
Then, by the choice of $\gamma$, there exists $h\in\alpha(\omega)$ such that
\begin{equation}\label{utv-e3}
\res{h}{H_2}=\gamma,
\end{equation}
and at the same time (see \fref{utv1-beta})
\begin{equation}\label{utv-e4}
\res{h}{H_1}\in K.
\end{equation}

From \fref{utv-e1} and the choice of $\omega$, $\omega'$ it follows that there is $h'\in\alpha(\omega')$ satisfying the equality
\begin{equation}\label{utv-e5}
\res{h'}{H_2}=\res{h}{H_2}.
\end{equation}
From \fref{utv-e3}, \fref{utv-e5} we have
\begin{equation}\label{utv-e6}
\res{h'}{H_2}=\gamma.
\end{equation}
In addition, \fref{utv-e4}, \fref{utv-e5}, and \fref{utv-e0} imply
\begin{equation}\label{utv-e7}
\res{h'}{H_1}\in K.
\end{equation}
From the inclusion of $h'\in\alpha(\omega')$, \fref{utv-e7} and the definition of $\beta$ (see \fref{utv1-beta}), we get $h'\in\beta(\omega')$, whence, taking into account \fref{utv-e6}, we have $\gamma\in\sres{\beta(\omega')}{H_2}$.
Then, due to the arbitrary choice of $\gamma$, we have $\sres{\beta(\omega)}{H_2}\subset\sres{\beta(\omega')}{H_2}$.
From here, in view of the symmetry of the occurrence of $\omega$, $\omega'$, the desired equality \fref{utv-e2} is extracted.
The proof is complete.

\bigskip

\begin{cor}\label{zam1}
In particular, if \fref{utv-e0} is true, the application of $\Ana\cdot{H_1}$ to $H_2$-non-anticipative m/f $\alpha$ does not violate this property.
\end{cor}

\proof
Indeed, let $H_1,H_2\in\naTc$ satisfy \fref{utv-e0}, $\alpha\in\AnaN{H_2}{}$, and $\omega\in\naOm$ is fixed.
Let us define $K$ and $\Upsilon$ from lemma \ref{utv1} as follows:
$$
K\mydef\bigcap_{\bar\omega\in\naOmo\omega{H_1}}\sres{\alpha(\bar\omega)}{H_1},\qquad\Upsilon\mydef\naOmo\omega{H_2}.
$$
Then, all conditions of lemma \ref{utv1} are satisfied and m/f $\beta$ specified in \fref{utv1-beta} satisfies the equalities (see \fref{def-Ana}) for all $\tilde\omega\in\Upsilon=\naOmo\omega{H_2}$:
\begin{multline*}
\beta(\tilde\omega)\mydef\{h\in\alpha(\tilde\omega)\mid\res h{H_1}\in K\}\mydef\{h\in\alpha(\tilde\omega)\mid\res h{H_1}\in\bigcap_{\bar\omega\in\naOmo{\omega}{H_1}}\sres{\alpha (\bar\omega)}{H_1}\}\\
=\{h\in\alpha(\tilde\omega)\mid\res h{H_1}\in\bigcap_{\bar\omega\in\naOmo{\tilde\omega}{H_1}}\sres{\alpha (\bar\omega)}{H_1}\}\mydef\Ana\alpha{H_1}(\tilde\omega).
\end{multline*}
Using the last relations and  \fref{utv-e2}, we obtain:
\beq\label{zam1-eq}
\sres{\Ana{\alpha}{H_1}(\omega'')}{H_2}=\sres{\beta(\omega'')}{H_2}=\sres{\beta(\omega')}{H_2}= \sres{\Ana{\alpha}{H_1}(\omega')}{H_2}\qquad\forall\omega',\omega''\in\Upsilon.
\eeq
Due to arbitrary choice of $\omega$, relations \fref{zam1-eq} are true for all $\omega',\omega''\in\naOmo\omega{H_2}$ and all $\omega\in\naOm$, i.\,e. we have inclusion $\Ana{\alpha}{H_1}\in\AnaN{H_2}{}$.
The proof is complete.

\bigskip

The following statement provides for a finite chain ${\cal H}$ and an arbitrary m/f a representation of its greatest ${\cal H}$-non-anticipative multiselector as a finite superposition of operators from ${\bf F}_{\cal H}$.
The provided construction inherits the features of the backward recurrent procedures \cite{NNK:143et,Pontr1967et,Pshenichny:1969et} on the one hand, and of the method of programmed iterations of A.\,G.\,Chentsov, on the other hand.

\begin{theorem}
\label{fix-to-couple-lem}
Let $\alpha\in\icP{\naZ}^\naOm$ and ${\cal H}=\{H_i\in\naTc\mid i\in\nint1k, k\in\NA\}$ be a finite chain: $(i\le j)\myeqv(H_i\subset H_j)$, $\forall i,j\in\nint1k$.
Then the equality is fulfilled:
\beq\label{fix-to-fin}
\Ana{\ldots\Ana\alpha{H_k}\ldots}{H_1}=\Ana\alpha{\cal H}.
\eeq
\end{theorem}

That is, the expression on the left gives (see \fref{Ana-max-H}) the greatest ${\cal H}$-non-anticipative multiselector of m/f $\alpha$.

\proof
1.
Denote $\phi\mydef\Ana{\ldots\Ana\alpha{H_k}\ldots}{H_1}$.
Then, successively applying the operators $\Ana\cdot{H_k}$,\ldots,$\Ana\cdot{H_1}$ to the m/f $\alpha$ and using the property \fref{Ana-cont} we arrive at the inequality
\beq\label{fi-in-al}
\phi\myLe\alpha;
\eeq
also applying successively the operators $\Ana\cdot{H_k}$,\ldots,$\Ana\cdot{H_1}$ to the inequality $\Ana\alpha{\cal H}\myLe\alpha$ taking into account isotonicity (see \fref{Ana-isot}) we get the ratio
\beq\label{ana-al-le-bet}
\Ana\alpha{\cal H}\myLe\phi.
\eeq
Therefore (see (ii) of lemma  \ref{lem-Ana-max}), the assertion will be proven if we establish the ${\cal H}$-non-anticipatory property of $\phi$.

2.
Let's show that for all $i\in\nint1k$ m/f $\phi$ holds the property \fref{fi-prop} of $H_i$-non-anticipatory:
\beq\label{fi-prop}
\sres{\phi(\omega)}{H_i}=\sres{\phi(\omega')}{H_i}, \qquad\forall \omega,\omega'\in\naOm,\ \omega'\in \naOmo{\omega}{H_i}.
\eeq

Case $i=k$.
M/f $\Ana{\alpha}{H_k}$ is $H_k$-non-anticipative by construction (see (i) of lemma \ref{lem-Ana-max}):
\beq\label{al-prop-k}
\sres{\Ana\alpha{H_k}(\omega)}{H_k}=\sres{\Ana\alpha{H_k}(\omega')}{H_k},\qquad\forall \omega,\omega'\in\naOm,\ \omega'\in\naOmo{\omega}{H_k}.
\eeq
Since ${\cal H}$ is a chain and hence $H_{k-1}\subset H_k$, taking into account the corollary \ref{zam1} from \fref{al-prop-k}, we get the property of $H_k$-non-anticipatory of m/f $\Ana{\Ana\alpha{H_k}}{H_{k-1}}$:
\beq\label{al-prop-k-1}
\sres{\Ana{\Ana\alpha{H_k}}{H_{k-1}}(\omega)}{H_k}=\sres{\Ana{\Ana\alpha{H_k}}{H_{k- 1}}(\omega')}{H_k},\qquad\forall \omega,\omega'\in\naOm,\ \omega'\in\naOmo{\omega}{H_k}.
\eeq
Due to relations \fref{al-prop-k-1} we can apply the reasoning to m/f $\Ana{\Ana\alpha{H_k}}{H_{k-1}}$ and operator $\Ana{\cdot}{H_{k-2}}$.
Continuing these arguments up to the application of the operator $\Ana{\cdot}{H_1}$ inclusively, we obtain the property of $H_k$-non-anticipatory for the m/f $\phi$, i.e., we prove statement \fref{fi-prop} for the case $i=k$.

Case $i\in\nint1{(k-1)}$.
Consider m/f $\psi\mydef\Ana{\ldots\Ana{\alpha}{H_k}\ldots}{H_i}$.
Since the operator $\Ana{\cdot}{H_i}$ was used last, by virtue of item (i) of lemma \ref{lem-Ana-max} $\psi$ is a $H_i$-non-anticipative m/f.
Further, repeating the arguments from the case $i=k$ for the operators $\Ana{\cdot}{H_j}$, $j\in\nint1{(i-1)}$, we conclude that these operators, applied to m/f $\psi$ when constructing m/f $\phi$, preserve the $H_i$-non-anticipatory property of m/f $\psi$.
Namely, for all $j\in\nint1{(i-1)}$ the equalities are hold true:
\beqnt
\sres{\Ana{\ldots\Ana{\psi(\omega)}{H_{i-1}}\ldots}{H_j}}{H_i}=\sres{\Ana{\ldots\Ana{\psi (\omega')}{H_{i-1}}\ldots}{H_j}}{H_i}, \qquad\forall \omega,\omega'\in\naOm,\ \omega'\in\naOmo{\omega}{H_i}.
\eeq
In particular, for $j=1$ we have \fref{fi-prop} for $i\in\nint1{(k-1)}$.
So, m/f $\phi$ is ${H_i}$-non-anticipative for all $i\in\nint1k$:
\beq\label{fi-nona}
\phi\in\AnaN{\cal H}{\alpha}.
\eeq
From the relations \fref{fi-in-al}, \fref{fi-nona} and item (ii) of lemma \ref{lem-Ana-max}, the relation $\phi\myLe\Ana\alpha{\cal H}$ follows.
Together with \fref{ana-al-le-bet} the relation gives us equality $\phi=\Ana\alpha{\cal H}$, i.e. the required equality \fref{fix-to-fin} is true.
The proof is complete.

\bigskip

\begin{remark}
Lemma \ref{utv1} and theorem \ref{fix-to-couple-lem} remain true in general case of sets  $\naT$ and $\naTc$, i.e. when $\naT$ is a set and  $\naTc$ is a chain in poset $(\icPp{\naT},\subset)$.
\end{remark}

\bigskip

Denote by $\stm\omega{\omega'}\in\icP\naTc$  the set of  the form (see \cite[(3.1),(3.6)]{GomSer-UDSU2021})
$$
\stm\omega{\omega'}\mydef\{H\in\naTc\mid\res\omega H=\res{\omega'}H \},\qquad\omega,\omega'\in\naOm,
$$
and define the mapping $\stmb\cdot\cdot\in\naTc^{\naOm\times\naOm}$ given by
$$
\stmb\omega{\omega'}\mydef\bigcup_{H\in\stm\omega{\omega'}}H,\qquad\omega,\omega'\in\naOm.
$$

The result of Corollary \ref{3netrudno} is close to the idea of a monotonicity of non-anticipatory m/fs (see \cite[Remark 2.8]{CARDALIAGUET-PLASKACZ-SIAM2000}).

\begin{cor}\label{3netrudno}
Let $\alpha\in\icP{\naZ}^\naOm$, $\naOm$ be a finite set, and
\beq\label{stnb-in-stm}
\stmb\omega{\omega'}\in\stm\omega{\omega'},\qquad\forall\omega,\omega'\in\naOm.
\eeq
Then the set $\naTc_\naOm$, $\naTc_\naOm\in\icP\naTc$, defined by $\naTc_\naOm\mydef\{\stmb\omega{\omega'}\,:\,\omega,\omega'\in\naOm\}$ is a finite chain, $\naTc_\naOm=\{\bar H_1\subset\ldots\subset\bar H_{k_{\naTc_\naOm}}\}$, and
\beq\label{tAA}
\top_{\naNo\alpha}=\Ana\alpha{\naTc_\naOm}=\Ana{\ldots\Ana\alpha{\bar H_{k_{\naTc_\naOm}}}\ldots}{\bar H_1}.
\eeq
\end{cor}

\proof
By definitions the inclusion $\naNo\alpha\subset\AnaN{\naTc_\naOm}\alpha$ is fulfilled (see \fref{nano-def}, sec. \ref{Ana}, \fref{A-nan-isot-H}).
On the other hand, suppose $\beta\in\AnaN{\naTc_\naOm}\alpha$ and $(\omega,\omega',H)\in\naOm\times\naOm\times\naTc$.
If $\res\omega H=\res{\omega'} H$, then $H\in\stm\omega{\omega'}$, $\stmb\omega{\omega'}\in\naTc_\naOm$, and $H\subset\stmb\omega{\omega'}$.
Then due to inclusion \fref{stnb-in-stm}, we have
$$
\res\omega{\stmb\omega{\omega'}}=\res{\omega'}{\stmb\omega{\omega'}}
$$
and by the property of $\naTc_\naOm$-non-anticipatory of $\beta$, the relation is fulfilled:
$$
\sres{\beta(\omega)}{\stmb\omega{\omega'}}=\sres{\beta(\omega')}{\stmb\omega{\omega'}}.
$$
Wherefrom the equality
$$
\sres{\beta(\omega)}{H}=\sres{\beta(\omega')}{H}
$$
immediately follows.
Since the choice of $(\omega,\omega',H)$ was arbitrary, we have the inclusion $\beta\in\naNo\alpha$, i.e. $\AnaN{\naTc_\naOm}\alpha\subset\naNo\alpha$.
Then, the equality $\naNo\alpha=\AnaN{\naTc_\naOm}\alpha$ is true.
By definition, this implies the first equality in \fref{tAA}.
The second equality in \fref{tAA} follows from finiteness of $\naTc_\naOm$ and theorem \ref{fix-to-couple-lem}.

The proof is complete.

\bigskip

Let us return to step-by-step finding of the selector under the conditions $(\alpha,\Delta)$.
Theorem \ref{fix-to-couple-lem} specifies a way to construct m/f $\Ana\alpha{{\cal H}_\Delta}$ --- the $\myLe$-greatest ${\cal H}_\Delta$-non-anticipative multiselector of $\alpha$.
If this multiselector turns out to be non-empty-valued, then, by virtue of theorem \ref{lem-aHfiz-eq-aHdOm}, the step-by-step procedure can be implemented by means of  $\Ana\alpha{{\cal H}_\Delta}$ for any disturbance $\omega\in\naOm$ (see the examples \ref{ex3}, \ref{ex4}).
At the same time, as in example \ref{ex3}, an ordinary non-anticipative multiselector of m/f $\alpha$ may be absent ($\naNO\alpha=\varnothing$).

In the case, when for some $\omega\in\naOm$ the value of the multiselector is empty ($\Ana\alpha{{\cal H}_\Delta}(\omega)=\varnothing$) due to $\myLe$-majority of $\Ana\alpha{{\cal H}_\Delta}$ in $\AnaN{{\cal H}_\Delta}\alpha$, we obtain the fact of unrealizability (theorem \ref{lem-aHfiz-eq-aHdOm}) of step by step procedure under the conditions $(\alpha,\Delta)$ and, as a consequence (see \fref{A-nan-isot-H}), under any other conditions $(\alpha,\Delta')$ where $\Delta\subset\Delta'$.

\section{Examples}
\label{EXAM}

\subsection{Example 1}
\label{ex1}

\begin{figure}[!ht]
\begin{center}
\hspace*{-0.5cm}
\input{pic1.tex}
\caption{Example 1}\label{pic1-ex1}
\end{center}
\end{figure}
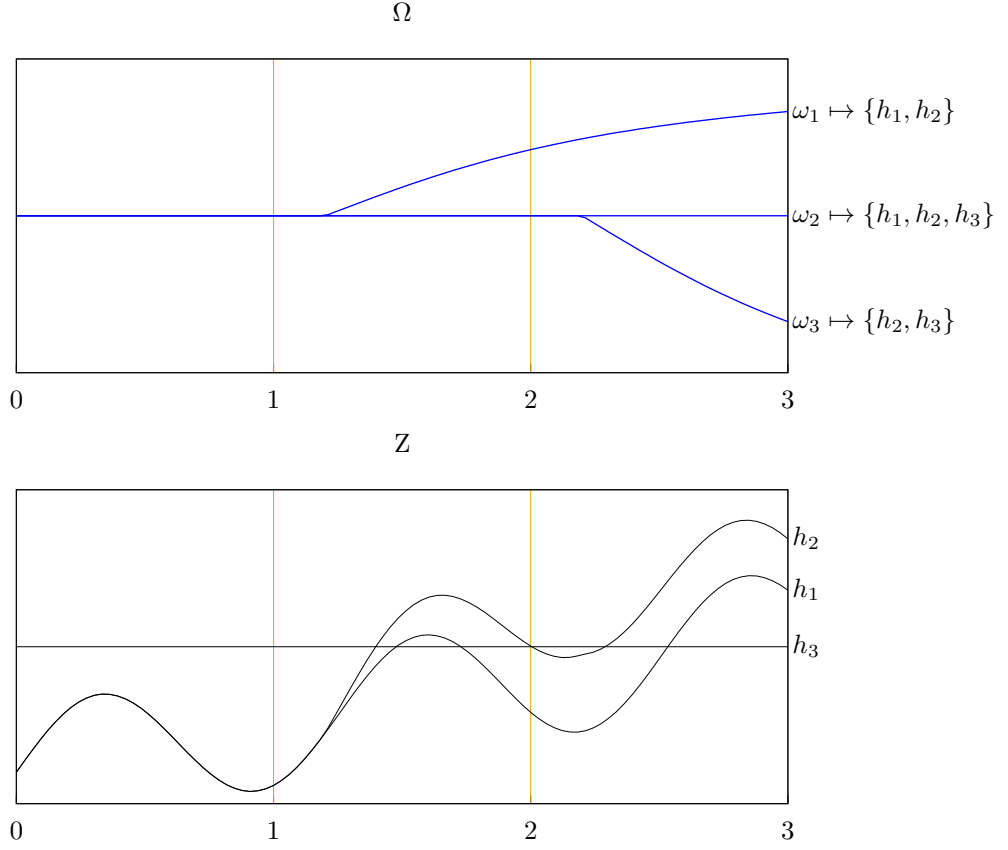

Let $\naT\mydef[0,3]$, $\naTc\mydef\{[0,\tau]\mid\tau\in\naT\}$ and $\naY=\naX=\RA$.
The sets $\naOm$ and $\naZ$ are shown in figure \ref{pic1-ex1}.
Let m/f $\beta$ be of the form
$$
\beta(\omega)=
\begin{cases}
\{h_1,h_2\},&\omega=\omega_1,\\
\{h_1,h_2,h_3\},&\omega=\omega_2,\\
\{h_2,h_3\},&\omega=\omega_3.
\end{cases}
$$
For this m/f, using definition \fref{def-Ana}, we get:
$$
\Ana{\beta}{[0,1]}(\omega)=
\begin{cases}
\{h_1,h_2\},&\omega=\omega_1,\\
\{h_1,h_2\},&\omega=\omega_2,\\
\{h_2\},&\omega=\omega_3,
\end{cases}
\qquad
\Ana{\beta}{[0,2]}(\omega)=
\begin{cases}
\{h_1,h_2\},&\omega=\omega_1,\\
\{h_2,h_3\},&\omega=\omega_2,\\
\{h_2,h_3\},&\omega=\omega_3.
\end{cases}
$$
It is clear that the inequalities $\Ana{\beta}{[0,1]}\not\myLe\Ana{\beta}{[0,2]}$, $\Ana{\beta}{[0,2 ]}\not\myLe\Ana{\beta}{[0,1]}$ are fulfilled.
That is, in general, the mapping $\naTc\ni H\mapsto\Ana{\beta}{H}\in{\icP\naZ}^\naOm$, considered as a mapping from poset $(\naTc,\subset)$ in poset $({\icP\naZ}^\naOm,\myLe)$ does not posses the isotonic property.

\subsection{Example 2}
\label{ex2}

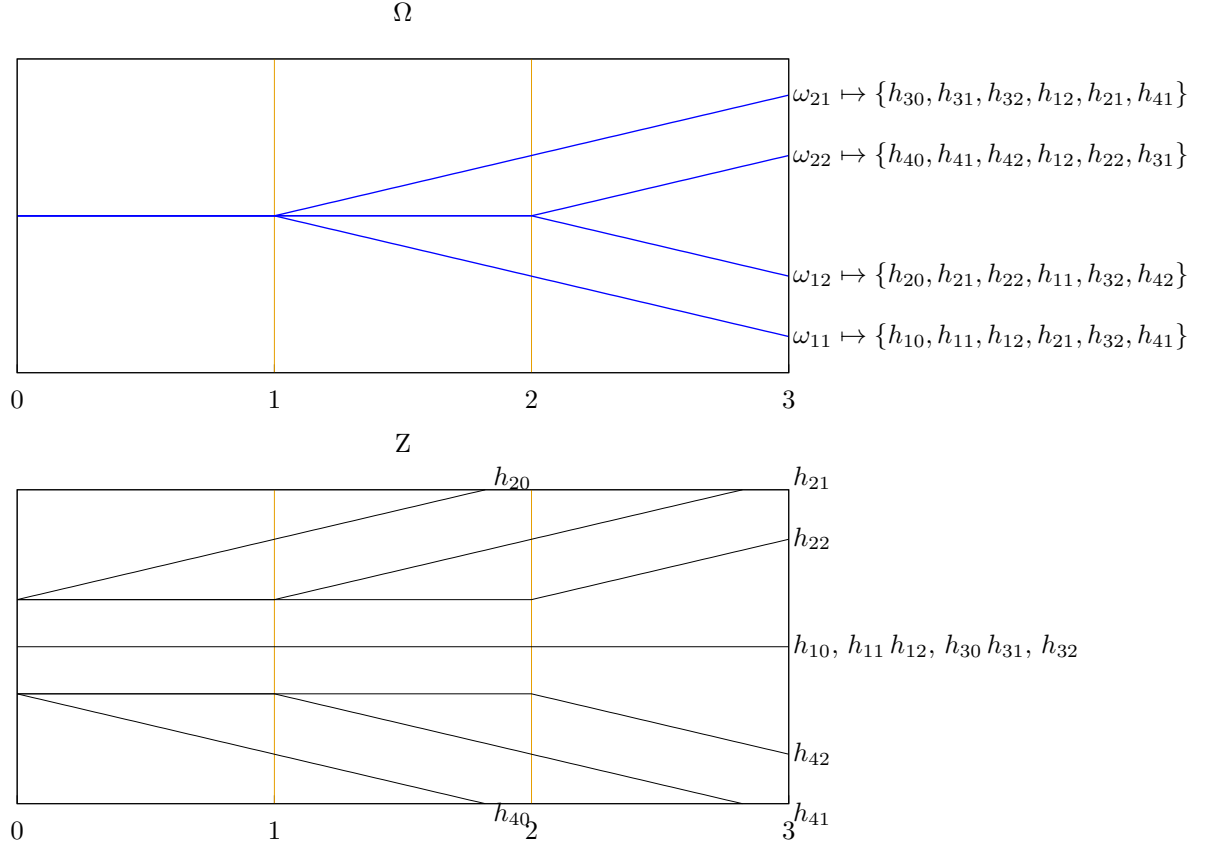
\begin{figure}[!ht]
\begin{center}
\hspace*{-0.5cm}\input{pic2.tex}
\caption{Example 2}\label{pic2-ex2}
\end{center}
\end{figure}

The example shows that the operators $\Ana{\cdot}{H}$, $H\in\naTc$ can be non-commutative.
Let $\naT\mydef[0,3]$, $\naTc\mydef\{[0,\tau]\mid\tau\in\naT\}$ and $\naY=\RA$, $\naX=X_1\times X_2=\RA^2$.
Denote
$$
\omega_{ij}(t)\mydef(-1)^i\max\{0,t-j\},\quad h_{ij}(t)\mydef a_i\left(1+\max\{0,t-j \}\right),\quad i\in\nint14,\ j\in\nint02,\ t\in[0,3],
$$
where $a_1\mydef(1,0)$, $a_2\mydef(0,1)$, $a_3\mydef(-1,0)$, $a_4\mydef(0,-1)$.
Put (see figure \ref{pic2-ex2}; for reasons of symmetry of $\naZ$, the projection onto the plane $\naT\times X_2$ is only shown)
$$
\naOm\mydef\{\omega_{ij}\mid i,j\in\nint12\},\qquad \naZ\mydef\{h_{ij}\mid i\in\nint14,\ j\in\nint02\}.
$$
Thus, all elements in $\naOm$ and $\naZ$ are distinct.
Moreover, if $k=i$, then
$$
\res{\omega_{kl}}{[0,\min\{l,j\}]}=\res{\omega_{ij}}{[0,\min\{l,j\}]}, \quad\res{h_{kl}}{[0,\min\{l,j\}]}=\res{h_{ij}}{[0,\min\{l,j\}]}.
$$
Let m/f $\alpha$ be given by
\begin{gather*}
   \alpha(\omega_{11})\mydef\{h_{10},h_{11},h_{12},h_{21},h_{32},h_{41}\},\quad
   \alpha( \omega_{21})\mydef\{h_{30},h_{31},h_{32},h_{12},h_{21},h_{41}\}, \\
   \alpha(\omega_{12})\mydef\{h_{20},h_{21},h_{22},h_{11},h_{32},h_{42}\},\quad
   \alpha(\omega_{22})\mydef\{h_{40},h_{41},h_{42},h_{12},h_{22},h_{31}\}.
\end{gather*}
Let $H_1\mydef[0,1]$.
Then $\naOm=\naOmo{\omega}{H_1}$ for all $\omega\in\naOm$;
hence, by definition \fref{def-Ana}, we have the equalities
\begin{gather*}
  \Ana\alpha{H_1}(\omega_{11})=\{h_{11},h_{12},h_{21},h_{32},h_{41}\}, \\
  \Ana\alpha{H_1}(\omega_{12})=\{h_{21},h_{22},h_{11},h_{32},h_{42}\}, \\
  \Ana\alpha{H_1}(\omega_{21})=\{h_{31},h_{32},h_{12},h_{21},h_{41}\}, \\
  \Ana\alpha{H_1}(\omega_{22})=\{h_{41},h_{42},h_{12},h_{22},h_{31}\}.
\end{gather*}
Let $H_2\mydef[0,2]$.
Then $\naOmo{\omega_{11}}{H_2}=\{\omega_{11}\}$, $\naOmo{\omega_{21}}{H_2}=\{\omega_{21}\}$, $\naOmo{\omega_{12}}{H_2}=\naOmo{\omega_{22}}{H_2}=\{\omega_{12},\omega_{22}\}$;
and by virtue of \fref{def-Ana} we get
\begin{eqnarray*}
  \Ana\alpha{H_2}(\omega_{11})&=&\{h_{10},h_{11},h_{12},h_{21},h_{32},h_{41}\} , \\
  \Ana\alpha{H_2}(\omega_{12})&=&\{h_{22},h_{42}\}, \\
  \Ana\alpha{H_2}(\omega_{21})&=&\{h_{30},h_{31},h_{32},h_{12},h_{21},h_{41}\} , \\
  \Ana\alpha{H_2}(\omega_{22})&=&\{h_{42},h_{22}\}.
\end{eqnarray*}
Then for m/f $\Ana{\Ana\alpha{H_1}}{H_2}$ we have
\begin{eqnarray*}
  \Ana{\Ana\alpha{H_1}}{H_2}(\omega_{11})&=&\{h_{11},h_{12},h_{21},h_{32},h_{41} \}, \\
  \Ana{\Ana\alpha{H_1}}{H_2}(\omega_{12})&=&\{h_{22},h_{42}\}, \\
  \Ana{\Ana\alpha{H_1}}{H_2}(\omega_{21})&=&\{h_{31},h_{32},h_{12},h_{21},h_{41} \}, \\
  \Ana{\Ana\alpha{H_1}}{H_2}(\omega_{22})&=&\{h_{42},h_{22}\}.
\end{eqnarray*}
Finally, for the m/f $\Ana{\Ana\alpha{H_2}}{H_1}$ we have the relations:
\begin{eqnarray*}
  \Ana{\Ana\alpha{H_2}}{H_1}(\omega_{11})&=&\{h_{21},h_{41}\}, \\
  \Ana{\Ana\alpha{H_2}}{H_1}(\omega_{12})&=&\{h_{22},h_{42}\}, \\
  \Ana{\Ana\alpha{H_2}}{H_1}(\omega_{21})&=&\{h_{21},h_{41}\}, \\
  \Ana{\Ana\alpha{H_2}}{H_1}(\omega_{22})&=&\{h_{42},h_{22}\}.
\end{eqnarray*}
It is easy to see that there is an inequality $\Ana{\Ana\alpha{H_2}}{H_1}\neq\Ana{\Ana\alpha{H_1}}{H_2}$, indicating that the operators $\Ana{\cdot}{H_1}$ and $\Ana{\cdot}{H_2}$ are not commutative.

Since the composition $\Ana{\Ana\cdot{H_2}}{H_1}$ corresponds to the order specified in theorem \ref{fix-to-couple-lem}, the result --- $\Ana{\Ana\alpha{H_2}}{H_1}$ --- presents an $\{H_1,H_2\}$-non-anticipative multiselector of $\alpha$.

\subsection{Example 3}
\label{ex3}

In this example, we consider an approaching game problem in which the m/f of optimal trajectories does not have any non-empty-valued and non-anticipative multiselector, that is, the problem is not solvable in the class of non-anticipative strategies built on usual (not relaxed) controls.
At the same time, the m/f of optimal trajectories has the property of feasibility for all partitions of the control interval;
i.\,e. the step-by-step procedure can be fulfilled for any partition and any disturbance.

\begin{figure}[!ht]
\begin{center}
\hspace*{-0.5cm}\input{pic3.tex}
\end{center}
\caption{Example 3}\label{pic3-ex3}
\end{figure}
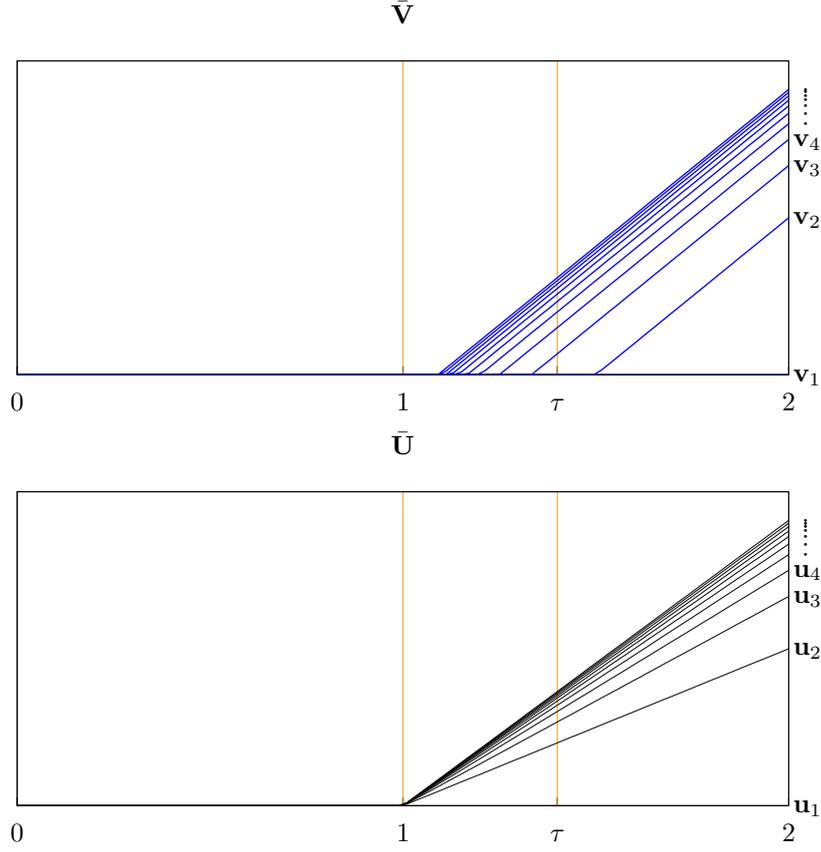

Let the trajectories of the controlled system be given by solutions of the following Cauchy problem:
\beq\label{sys5}
\begin{cases}
\dot x(t)=u(t)-v(t),& t\in\naT\mydef[0,2],\\
x(0)=0\in\RA,&u\in\UB, v\in\VB,
\end{cases}
\eeq
\beq\label{UV-def5}
\UB\mydef\{u_i\mid i\in\NA\},\qquad \VB\mydef\{v_i\mid i\in\NA\},
\eeq
\beq\label{UV-def51}
u_i(t)\mydef
\begin{cases}
0,&t\in[0,1],\\
1-1/i,&t\in(1,2],
\end{cases}
\qquad
v_i(t)\mydef
\begin{cases}
0,&t\in[0,1+1/i],\\
1,&t\in(1+1/i,2],
\end{cases}
\qquad i\in\NA.
\eeq

For system \fref{sys5}--\fref{UV-def51}, consider the problem of meeting its trajectories $x(\cdot)$ with the set $M\mydef\{(t,x)\in\naT\times\RA\mid (t=2)\myand (x\ge0)\}$ by choosing of control $u(\cdot)\in\UB$ for all possible disturbances $v(\cdot)\in\VB$.
Denote by $x(\cdot,u,v)$, where $x(\cdot,u,v)\in C(\naT,\RA)$, the solution of the Cauchy problem \fref{sys5}--\fref{UV-def51} where a control $u\in\UB$ and a disturbance $v\in\VB$ are given.
Denote $\UBB\mydef\{{\bf u}_i\mid i\in\NA\}$, $\VBB\mydef\{{\bf v}_i\mid i\in\NA\}$ (see figure\,\ref{pic3-ex3})
$$
{\bf u}_i(t)\mydef\int_{0}^{t}u_i(s)ds,\quad{\bf v}_i(t)\mydef\int_{0}^{t}v_i(s)ds,\qquad i\in\NA,\ t\in\naT.
$$
It is easy to verify that for $x_{ij}(\cdot)\mydef x(\cdot,u_i,v_j)$ the equalities (see \fref{UV-def5}, \fref{UV-def51})
$$
x_{ij}(t)={\bf u}_i(t)-{\bf v}_j(t),\quad{\bf u}_i(2)={\bf v}_i(2)=1-\frac1i,\qquad i,j\in\NA, t\in\naT,
$$
and, consequently, the equalities
\beq\label{x-ij}
x_{ij}(2)=\frac1j-\frac1i,\qquad i,j\in\NA
\eeq
are fulfilled.
Hence, given a disturbance $v_j(\cdot)\in\VB$ and a control $u_i(\cdot)$, the meeting criterion can be written as the inequality
$$
i\ge j.
$$
Therefore, for any disturbance $v_j(\cdot)\in\VB$, denoting by $\alpha(v_j(\cdot))$ the set of all controls in $\UB$ that solve the meeting problem, we can write
\beq\label{al5-def}
\alpha(v_j)=\{u_j,u_{j+1},\ldots\}=\{u_i\in\UB\mid i\ge j\}.
\eeq
That is, in \fref{al5-def} we have the m/f of optimal answers.
The family of subsets $\naTc$, implementing "flow of time", has its usual form: $\naTc\mydef\{[0, \tau]\mid \tau\in\naT\}$.

We choose a finite partition $\Delta\subset\naT$ and show that conditions $(\alpha,\Delta)$ are feasible.
To this end, for ${\cal H}_\Delta\in\naTc$ we construct using \fref{fix-to-fin} ${\cal H}_\Delta$-non-anticipative multiselector of $\alpha$ and verify that its values be non-empty.
Thus, by virtue of \fref{aHfiz-eq-aHdOm}, the feasibility of step-by-step procedure under conditions $(\alpha,\Delta)$ will be proven for arbitrary finite partition $\Delta$.

1.
Denote $H_\tau\mydef[0,\tau]$.
Then
$$
\naOmo{v}{H_\tau}=\VB, \qquad \tau\in[0,1],\ v\in\VB;
$$
$$
\naOmo{v_j}{H_\tau}=
\begin{cases}
\{v_i,\ldots,v_1\},& j\le i,\\
\{v_j\},& j>i,
\end{cases}
\qquad \tau\in\left(1+\frac1{i+1},1+\frac1i\right],\quad i, j\in\NA.
$$

Hence, by direct calculation (see \fref{def-Ana}, \fref{al5-def}), we obtain the representation of m/f $\Ana{\alpha}{H_\tau}$:
\begin{gather}
\Ana{\alpha}{H_\tau}=\alpha, \qquad \tau\in[0,1],\label{ana-al0-ex5}\\
\Ana{\alpha}{H_\tau}(v_j)=
\begin{cases}
\alpha(v_i),&j\le i,\\
\alpha(v_j),&j>i,
\end{cases}
\qquad \tau\in\left(1+\frac1{i+1},1+\frac1i\right],\quad i, j\in\NA.\label{ana-al1-ex5}
\end{gather}

From the equalities \fref{ana-al0-ex5} and \fref{ana-al1-ex5} the relations follow
\beq\label{ana-al2-ex5}
\nado{\Ana{\alpha}{H_\tau}}=\naOm,\qquad\forall\tau\in\naT.
\eeq
Moreover, it is not difficult to check that the representations \fref{ana-al0-ex5} and \fref{ana-al1-ex5} hold for a wider set of m/f, namely, for any $\beta\myLe\alpha$ take place formulas:
\beq\label{ana-be0-ex5}
\Ana{\beta}{H_\tau}=\beta, \qquad \tau\in[0,1],
\eeq
and for any $\beta$ of the form $\beta\mydef\Ana{\alpha}{H_\xi}$, $\xi\in\naT$, the equalities are true:
\beq\label{ana-be1-ex5}
\Ana{\beta}{H_\tau}(v_j)=
\begin{cases}
\beta(v_i),&i\ge j,\\
\beta(v_j),&i<j,
\end{cases}
\quad\ \tau\in\left(1+\frac1{i+1},1+\frac1i\right],\quad i,j\in\NA.
\eeq

2. Using presentations \fref{ana-al0-ex5}--\fref{ana-be1-ex5} in accordance with the \fref{fix-to-fin}, we finally obtain the equality
\beq\label{Ana-com5}
\Ana{\alpha}{{\cal H}_\Delta}=\Ana{\alpha}{[0,\min\{\Delta\cap(1,2]\}]}.
\eeq
Sines $\Delta$ is finite and $\Delta\cap(1,2]\neq\varnothing$ we have $\Ana{\alpha}{{\cal H}_\Delta}=\Ana{\alpha}{H_{\bar\tau}}$ for same $H_{\bar\tau}\in{\cal H}_\Delta$.
Then, in view of \fref{ana-al2-ex5}, we receive non-emptiness of values of m/f $\Ana{\alpha}{{\cal H}_\Delta}$.
That is, non-emptiness of the set $\AnaNO{{\cal H}_\Delta}\alpha$ of all non-empty-valued and ${{\cal H}_\Delta}$-non-anticipative multiselectors of $\alpha$.
Hence, taking into account theorem \ref{lem-aHfiz-eq-aHdOm}, we conclude that conditions $(\alpha, \Delta)$ are feasible.

3.
Let us show that the problem has no solution in the class of non-anticipative strategies (quasi-strategies), that is, m/f $\alpha$ \fref{al5-def} does not have non-empty-valued non-anticipative multiselectors: $\naNO\alpha=\varnothing$.

Let's say the contrary --- m/f $\beta\in\naNO\alpha$ was found.
Then \fref{A-nan-isot-H} implies that $\AnaNO\naTc\alpha\subset\AnaNO{\cal H}\alpha$, where ${\cal H}\mydef\{H_k\mid k\in\NA\}$ and $H_k\mydef[0,1+1/k]$.
Let's use the implication \fref{Ana-max-est}; by assumption we have $\beta\in\AnaNO\naTc\alpha\subset\AnaNO{\cal H}\alpha$, hence $\beta\in\bigwedge_{k\in\NA}\Ana{\alpha}{ H_k}$.
Then, taking into account \fref{ana-al1-ex5} and \fref{al5-def}, for each $j\in\NA$ we obtain
$$
\beta(v_j(\cdot))\subset\bigcap_{k\in\NA}\Ana{\alpha}{H_k}(v_j(\cdot))=\bigcap_{k>j}\alpha(v_k(\cdot))=\varnothing.
$$
So, $\nado{\beta}=\varnothing$ and relation $\beta\in\naNO\alpha$ is impossible.
Then, we have $\naNO\alpha=\varnothing$.

\subsection{Example 4}
\label{ex4}

Consider an example from \cite[\S5]{GomSer-UDSU2021}.
Here, we have the case where the construction of partially non-anticipative multiselector provides a classic non-anticipative multiselector.

As in previous example, the controlled system is given by solutions of the following Cauchy problem:
\beqnt
         \dot x(t)=u(t) + v(t), \quad t \in\naT \mydef [0, 3], \quad x(0) = 0,
\eeq
where $x(t) \in \mathbb{R}$, the control $u$ and disturbance $v$ are Borel measurable functions subject to the instantaneous constraints
\beqnt
         u(t) \in P \mydef [- 1, 1], \quad v(t) \in Q \mydef \{- 1, 0, 1\}, \quad t \in\naT.
\eeq
Denote by $x(\cdot,u,v)$, where $x(\cdot,u,v)\in C(\naT,\RA)$, the solution of the Cauchy problem where a control $u\in\UB$ and a disturbance $v\in\VB$ are given.

Suppose that the set $\VB$ of admissible disturbances consists of two functions, $\VB\mydef\{v_1, v_2\}$, where:
     \begin{gather*}
         v_1(t)
         \mydef \begin{cases}
             0, & \mbox{if } t \in [0, 1] \cup (2, 3], \\
             1, & \mbox{if } t \in (1, 2],
         \end{cases}
         \quad v_2(t)
         \mydef \begin{cases}
             0, & \mbox{if } t \in [0, 1], \\
             - 1, & \mbox{if } t \in (1, 3].
         \end{cases}
     \end{gather*}
The example meet the conditions of Corollary \ref{3netrudno}, in particular, the set of disturbances is finite.
So, due to the corollary, the class of non-anticipative strategies in this problem is limited by the set of all $[0,1]$-non-anticipative m/f  (from $\icPp{\UB}^\VB$);
here $\UB$ is the set of all possible realizations of control.

The control minimizes the following cost functional:
\beqnt
         J(u,v)\mydef - |x(3; u, v)|,\quad u \in \UB, \quad v\in \VB.
\eeq

Namely, the goal of control side is to minimize the guaranteed result in the class of non-anticipative strategies from $\icPp{\UB}^\VB$.

Since we are considering a guaranteed state of the problem, denote by $\rho$ (yet unknown) optimal (minimal) guaranteed result.
Let us compose the m/f $\alpha_\rho\in\icP{\UB}^\VB$, which describes the $\rho$-optimal control responses to the realized disturbance.

Taking into account the obvious inequality $\rho\le0$, we write the values of $\alpha_\rho$ depending on $\rho$:
\beq\label{ex40}
   \alpha_\rho(v)=\left\{u\in\UB\mid\left|\int_0^3u(s)+v(s)\dup s\right|\ge-\rho\right\},\qquad v\in\VB.
\eeq
Then we get the equalities
\beq\label{ex41}
   \alpha_\rho(v_1)=\left\{u\in\UB\mid\left(\int_0^3u(s)\dup s\le\rho-1\right)\vee\left(-\rho-1\le\int_0^3u(s)\dup s\right)\right\},
\eeq
\beq\label{ex42}
   \alpha_\rho(v_2)=\left\{u\in\UB\mid\left(\int_0^3u(s)\dup s\le\rho+2\right)\vee\left(-\rho+2\le\int_0^3u(s)\dup s\right)\right\}.
\eeq

It is clear that the set of values of the parameter $\rho$ for which $\alpha_{\rho}$ has a non-empty-valued and non-anticipative multiselector is bounded.
That is, the set has infimum $\bar\rho$, which is the optimal result of the control side in the class of non-anticipative strategies.
So, keeping in mind Corollary  \ref{3netrudno}, we can write:
$$
\bar\rho=\min\left\{\rho\in\RA\mid\nado{\top_{\naNo{\alpha_{\rho}}}}=\VB \right\}=\min\left\{\rho\in\RA\mid\nado{\Ana{\alpha_{\rho}}{\{[0,1]\}}}=\VB \right\}.
$$
Let us calculate the m/f $\Ana{\alpha_{\rho}}{\{[0,1]\}}$:
$$
\Ana{\alpha_{\rho}}{\{[0,1]\}}(v)\mydef\left\{u\in\alpha_\rho(v)\mid \res{u}{[0,1 ]}\in\bigcap_{v=v_1,v_2}\sres{\alpha_\rho(v)}{[0,1]}\right\},\qquad v\in\{v_1,v_2\}.
$$
Note that the condition $\nado{\Ana{\alpha_{\rho}}{\{[0,1]\}}}=\VB$ is equivalent to the inequality $\bigcap_{v=v_1,v_2}\sres{\alpha_\rho(v)}{[0,1]}\neq\varnothing$.
The last inequality is true if and only if (see \fref{ex41}, \fref{ex42}) the following two relations are fulfilled:
\begin{gather*}
\left (\rho-1-\min_{u\in\UB}\int_1^3u(s)\dup s\ge\int_0^1u(s)\dup s\right)\vee\left(\int_0^1u(s)\dup s\ge-\rho-1-\max_{u\in\UB}\int_1^3u(s)\dup s\right), \\
  \left (\rho+2-\min_{u\in\UB}\int_1^3u(s)\dup s\ge\int_0^1u(s)\dup s\right)\vee\left(\int_0^1u(s)\dup s\ge-\rho+2-\max_{u\in\UB}\int_1^3u(s)\dup s\right),
\end{gather*}
that is, the relations below are true:
\begin{gather*}
   \left( \rho+1\ge\int_0^1u(s)\dup s\right)\vee\left(\int_0^1u(s)\dup s\ge-\rho-3\right), \\
  \left ( \rho+4\ge\int_0^1u(s)\dup s\right)\vee\left(\int_0^1u(s)\dup s\ge-\rho\right).
\end{gather*}
Hence, we come to the equality
$$
\bar\rho=\min\left\{\rho\in\RA\mid\left\{u\in\UB\mid \rho+4\ge\int_0^1u(s)\dup s\ge-\rho-3\right\}\neq\varnothing\right\}.
$$
The latter condition implies that $\rho+4\ge-\rho-3$ or $\rho\ge-3.5$.

So, for the value $\bar\rho$ of optimal guarantee in the class of non-anticipative strategies, we have the equality $\bar\rho=-3.5$.
In this case, the greatest $[0,1]$-non-anticipative multiselector of $\alpha_{\bar\rho}$ is non-empty-valued and takes the form
  \begin{gather*}
         \Ana{\alpha_{\bar\rho}}{[0,1]}(v_1)(t)
         \mydef \begin{cases}
             u(t),\ u\in\UB_{([0,1],0.5)}, &t \in [0, 1], \\
             1=\argmax_{w\in P} w, &t \in (1, 3],
         \end{cases}\\
         \Ana{\alpha_{\bar\rho}}{[0,1]}(v_2)(t)
         \mydef \begin{cases}
             u(t),\ u\in\UB_{([0,1],0.5)}, &t \in [0, 1], \\
             - 1=\argmin_{w\in P}w, & t \in (1, 3],
         \end{cases}
     \end{gather*}
where
\begin{multline*}
\UB_{([0,1],0.5)}\mydef\Big\{u\in\UB\mid\res u{[0,1]}\in\bigcap_{v\in\{v_1,v_2\}}\sres{\alpha_{\bar\rho}(v)}{[0,1]}\Big\}\\
=\Big\{u\in\UB\mid\int_0^1u(s)\dup s=\bar\rho+4=-\bar\rho-3=0.5\Big\}.
\end{multline*}

Unlike Example \ref{ex3}, the optimal result here is achieved in the class of non-anticipative strategies, and can be explicitly described due to finiteness of the disturbances set.

\section{Conclusion}
\label{conc}
Since the topic of the paper is at the initial stage of development, many details and even essential issues remain unexplored.
Here, we make some remarks on them.
Expression \fref{fix-to-fin}, in view of the maximality properties of \fref{Ana-max} and \fref{Ana-max-H}, gives some hope for new results on existence of a m/f non-anticipative multiselector.
Concerning applications, an interesting question arises about the convergence of guaranteed results for ${\cal H}_\Delta$-non-anticipatory multiselectors of a m/f $\alpha$ to the guaranteed result of the non-anticipative multiselector of $\alpha$ as the step of the partition $\Delta$ tends to zero.
Another question is the implementation of the proposed constructions in a solution of dynamic optimization problem.
Namely, we need a systematic approach to resolving of the mathematical programming and the parametric optimization problems that arise when "calculating" partially non-anticipatory multiselectors (see the last example).

\bibliography{C:/Dropbox/CURRENT-WORK/ref/mybib,C:/Dropbox/CURRENT-WORK/ref/allbib}

\end{document}

%% file: pic1.tex
\begingroup
  \makeatletter
  \providecommand\color[2][]{%
    \GenericError{(gnuplot) \space\space\space\@spaces}{%
      Package color not loaded in conjunction with
      terminal option `colourtext'%
    }{See the gnuplot documentation for explanation.%
    }{Either use 'blacktext' in gnuplot or load the package
      color.sty in LaTeX.}%
    \renewcommand\color[2][]{}%
  }%
  \providecommand\includegraphics[2][]{%
    \GenericError{(gnuplot) \space\space\space\@spaces}{%
      Package graphicx or graphics not loaded%
    }{See the gnuplot documentation for explanation.%
    }{The gnuplot epslatex terminal needs graphicx.sty or graphics.sty.}%
    \renewcommand\includegraphics[2][]{}%
  }%
  \providecommand\rotatebox[2]{#2}%
  \@ifundefined{ifGPcolor}{%
    \newif\ifGPcolor
    \GPcolorfalse
  }{}%
  \@ifundefined{ifGPblacktext}{%
    \newif\ifGPblacktext
    \GPblacktexttrue
  }{}%
  \let\gplgaddtomacro\g@addto@macro
  \gdef\gplbacktext{}%
  \gdef\gplfronttext{}%
  \makeatother
  \ifGPblacktext
    \def\colorrgb#1{}%
    \def\colorgray#1{}%
  \else
    \ifGPcolor
      \def\colorrgb#1{\color[rgb]{#1}}%
      \def\colorgray#1{\color[gray]{#1}}%
      \expandafter\def\csname LTw\endcsname{\color{white}}%
      \expandafter\def\csname LTb\endcsname{\color{black}}%
      \expandafter\def\csname LTa\endcsname{\color{black}}%
      \expandafter\def\csname LT0\endcsname{\color[rgb]{1,0,0}}%
      \expandafter\def\csname LT1\endcsname{\color[rgb]{0,1,0}}%
      \expandafter\def\csname LT2\endcsname{\color[rgb]{0,0,1}}%
      \expandafter\def\csname LT3\endcsname{\color[rgb]{1,0,1}}%
      \expandafter\def\csname LT4\endcsname{\color[rgb]{0,1,1}}%
      \expandafter\def\csname LT5\endcsname{\color[rgb]{1,1,0}}%
      \expandafter\def\csname LT6\endcsname{\color[rgb]{0,0,0}}%
      \expandafter\def\csname LT7\endcsname{\color[rgb]{1,0.3,0}}%
      \expandafter\def\csname LT8\endcsname{\color[rgb]{0.5,0.5,0.5}}%
    \else
      \def\colorrgb#1{\color{black}}%
      \def\colorgray#1{\color[gray]{#1}}%
      \expandafter\def\csname LTw\endcsname{\color{white}}%
      \expandafter\def\csname LTb\endcsname{\color{black}}%
      \expandafter\def\csname LTa\endcsname{\color{black}}%
      \expandafter\def\csname LT0\endcsname{\color{black}}%
      \expandafter\def\csname LT1\endcsname{\color{black}}%
      \expandafter\def\csname LT2\endcsname{\color{black}}%
      \expandafter\def\csname LT3\endcsname{\color{black}}%
      \expandafter\def\csname LT4\endcsname{\color{black}}%
      \expandafter\def\csname LT5\endcsname{\color{black}}%
      \expandafter\def\csname LT6\endcsname{\color{black}}%
      \expandafter\def\csname LT7\endcsname{\color{black}}%
      \expandafter\def\csname LT8\endcsname{\color{black}}%
    \fi
  \fi
    \setlength{\unitlength}{0.0500bp}%
    \ifx\gptboxheight\undefined%
      \newlength{\gptboxheight}%
      \newlength{\gptboxwidth}%
      \newsavebox{\gptboxtext}%
    \fi%
    \setlength{\fboxrule}{0.5pt}%
    \setlength{\fboxsep}{1pt}%
\begin{picture}(6480.00,6480.00)%
    \gplgaddtomacro\gplbacktext{%
      \csname LTb\endcsname
      \put(330,3460){\makebox(0,0){\strut{}$0$}}%
      \put(2248,3460){\makebox(0,0){\strut{}$1$}}%
      \put(4165,3460){\makebox(0,0){\strut{}$2$}}%
      \put(6083,3460){\makebox(0,0){\strut{}$3$}}%
      \put(6083,5644){\makebox(0,0)[l]{\strut{}$\,\omega_1\mapsto\{h_1,h_2\}$}}%
      \put(6083,4860){\makebox(0,0)[l]{\strut{}$\,\omega_2\mapsto\{h_1,h_2,h_3\}$}}%
      \put(6083,4064){\makebox(0,0)[l]{\strut{}$\,\omega_3\mapsto\{h_2,h_3\}$}}%
    }%
    \gplgaddtomacro\gplfronttext{%
      \csname LTb\endcsname
      \put(3206,6369){\makebox(0,0){\strut{}$\naOm$}}%
    }%
    \gplgaddtomacro\gplbacktext{%
      \csname LTb\endcsname
      \put(330,220){\makebox(0,0){\strut{}$0$}}%
      \put(2248,220){\makebox(0,0){\strut{}$1$}}%
      \put(4165,220){\makebox(0,0){\strut{}$2$}}%
      \put(6083,220){\makebox(0,0){\strut{}$3$}}%
      \put(6083,2045){\makebox(0,0)[l]{\strut{}${\,h_1}$}}%
      \put(6083,2431){\makebox(0,0)[l]{\strut{}${\,h_2}$}}%
      \put(6083,1620){\makebox(0,0)[l]{\strut{}${\,h_3}$}}%
    }%
    \gplgaddtomacro\gplfronttext{%
      \csname LTb\endcsname
      \put(3206,3130){\makebox(0,0){\strut{}$\naZ$}}%
    }%
    \gplbacktext
    \put(0,0){\includegraphics{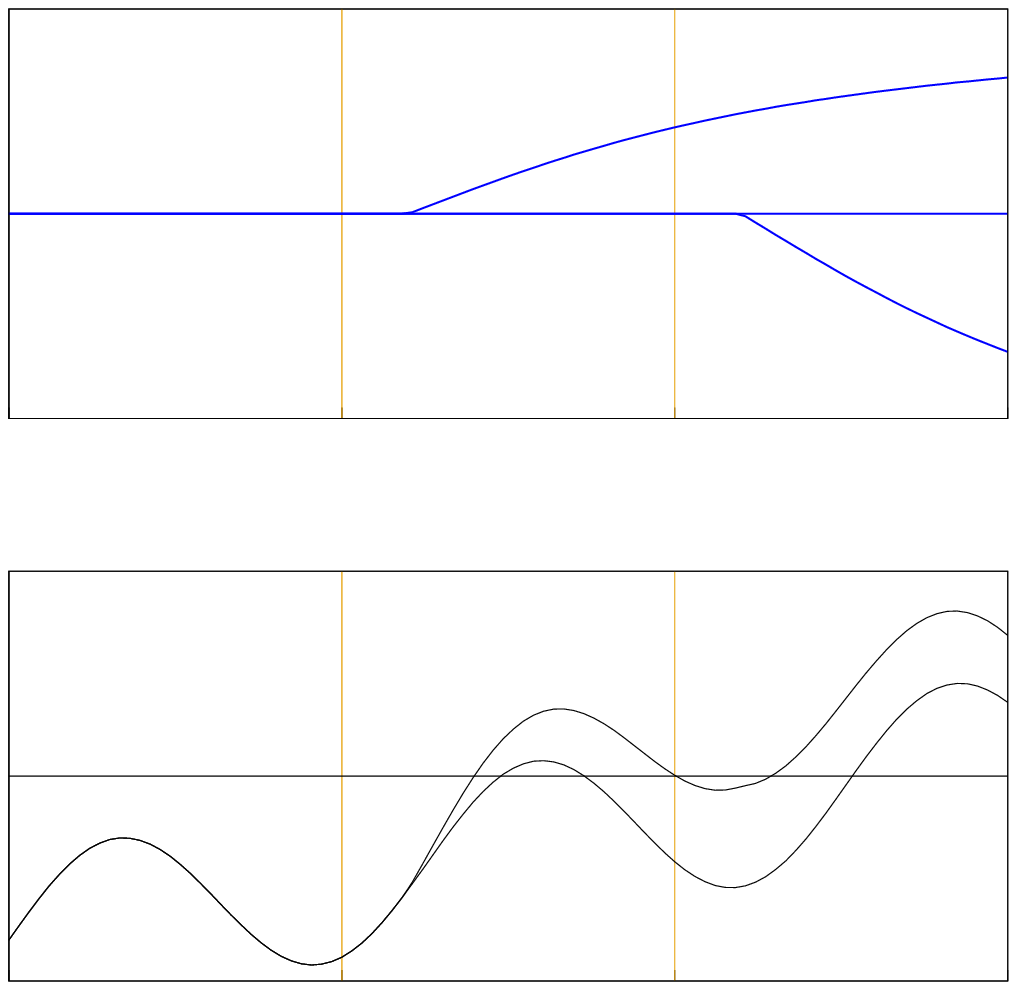}}%
    \gplfronttext
  \end{picture}%
\endgroup

%% file: pic2.tex
\begingroup
  \makeatletter
  \providecommand\color[2][]{%
    \GenericError{(gnuplot) \space\space\space\@spaces}{%
      Package color not loaded in conjunction with
      terminal option `colourtext'%
    }{See the gnuplot documentation for explanation.%
    }{Either use 'blacktext' in gnuplot or load the package
      color.sty in LaTeX.}%
    \renewcommand\color[2][]{}%
  }%
  \providecommand\includegraphics[2][]{%
    \GenericError{(gnuplot) \space\space\space\@spaces}{%
      Package graphicx or graphics not loaded%
    }{See the gnuplot documentation for explanation.%
    }{The gnuplot epslatex terminal needs graphicx.sty or graphics.sty.}%
    \renewcommand\includegraphics[2][]{}%
  }%
  \providecommand\rotatebox[2]{#2}%
  \@ifundefined{ifGPcolor}{%
    \newif\ifGPcolor
    \GPcolorfalse
  }{}%
  \@ifundefined{ifGPblacktext}{%
    \newif\ifGPblacktext
    \GPblacktexttrue
  }{}%
  \let\gplgaddtomacro\g@addto@macro
  \gdef\gplbacktext{}%
  \gdef\gplfronttext{}%
  \makeatother
  \ifGPblacktext
    \def\colorrgb#1{}%
    \def\colorgray#1{}%
  \else
    \ifGPcolor
      \def\colorrgb#1{\color[rgb]{#1}}%
      \def\colorgray#1{\color[gray]{#1}}%
      \expandafter\def\csname LTw\endcsname{\color{white}}%
      \expandafter\def\csname LTb\endcsname{\color{black}}%
      \expandafter\def\csname LTa\endcsname{\color{black}}%
      \expandafter\def\csname LT0\endcsname{\color[rgb]{1,0,0}}%
      \expandafter\def\csname LT1\endcsname{\color[rgb]{0,1,0}}%
      \expandafter\def\csname LT2\endcsname{\color[rgb]{0,0,1}}%
      \expandafter\def\csname LT3\endcsname{\color[rgb]{1,0,1}}%
      \expandafter\def\csname LT4\endcsname{\color[rgb]{0,1,1}}%
      \expandafter\def\csname LT5\endcsname{\color[rgb]{1,1,0}}%
      \expandafter\def\csname LT6\endcsname{\color[rgb]{0,0,0}}%
      \expandafter\def\csname LT7\endcsname{\color[rgb]{1,0.3,0}}%
      \expandafter\def\csname LT8\endcsname{\color[rgb]{0.5,0.5,0.5}}%
    \else
      \def\colorrgb#1{\color{black}}%
      \def\colorgray#1{\color[gray]{#1}}%
      \expandafter\def\csname LTw\endcsname{\color{white}}%
      \expandafter\def\csname LTb\endcsname{\color{black}}%
      \expandafter\def\csname LTa\endcsname{\color{black}}%
      \expandafter\def\csname LT0\endcsname{\color{black}}%
      \expandafter\def\csname LT1\endcsname{\color{black}}%
      \expandafter\def\csname LT2\endcsname{\color{black}}%
      \expandafter\def\csname LT3\endcsname{\color{black}}%
      \expandafter\def\csname LT4\endcsname{\color{black}}%
      \expandafter\def\csname LT5\endcsname{\color{black}}%
      \expandafter\def\csname LT6\endcsname{\color{black}}%
      \expandafter\def\csname LT7\endcsname{\color{black}}%
      \expandafter\def\csname LT8\endcsname{\color{black}}%
    \fi
  \fi
    \setlength{\unitlength}{0.0500bp}%
    \ifx\gptboxheight\undefined%
      \newlength{\gptboxheight}%
      \newlength{\gptboxwidth}%
      \newsavebox{\gptboxtext}%
    \fi%
    \setlength{\fboxrule}{0.5pt}%
    \setlength{\fboxsep}{1pt}%
\begin{picture}(6480.00,6480.00)%
    \gplgaddtomacro\gplbacktext{%
      \csname LTb\endcsname
      \put(330,3460){\makebox(0,0){\strut{}$0$}}%
      \put(2248,3460){\makebox(0,0){\strut{}$1$}}%
      \put(4165,3460){\makebox(0,0){\strut{}$2$}}%
      \put(6083,3460){\makebox(0,0){\strut{}$3$}}%
      \put(6083,3952){\makebox(0,0)[l]{\strut{}$\,\omega_{11}\mapsto\{h_{10},h_{11},h_{12},h_{21},h_{32},h_{41}\}$}}%
      \put(6083,4406){\makebox(0,0)[l]{\strut{}$\,\omega_{12}\mapsto\{h_{20},h_{21},h_{22},h_{11},h_{32},h_{42}\}$}}%
      \put(6083,5767){\makebox(0,0)[l]{\strut{}$\,\omega_{21}\mapsto\{h_{30},h_{31},h_{32},h_{12},h_{21},h_{41}\}$}}%
      \put(6083,5313){\makebox(0,0)[l]{\strut{}$\,\omega_{22}\mapsto\{h_{40},h_{41},h_{42},h_{12},h_{22},h_{31}\}$}}%
    }%
    \gplgaddtomacro\gplfronttext{%
      \csname LTb\endcsname
      \put(3206,6369){\makebox(0,0){\strut{}$\naOm$}}%
    }%
    \gplgaddtomacro\gplbacktext{%
      \csname LTb\endcsname
      \put(330,220){\makebox(0,0){\strut{}$0$}}%
      \put(2248,220){\makebox(0,0){\strut{}$1$}}%
      \put(4165,220){\makebox(0,0){\strut{}$2$}}%
      \put(6083,220){\makebox(0,0){\strut{}$3$}}%
      \put(6083,1620){\makebox(0,0)[l]{\strut{}${\,h_{10},\,h_{11}\,h_{12},\,h_{30}\,h_{31},\,h_{32}}$}}%
      \put(4165,2882){\makebox(0,0)[r]{\strut{}${\!\!h_{20}}$}}%
      \put(4165,358){\makebox(0,0)[r]{\strut{}${\!\!h_{40}}$}}%
      \put(6083,2882){\makebox(0,0)[l]{\strut{}${\,h_{21}}$}}%
      \put(6083,358){\makebox(0,0)[l]{\strut{}${\,h_{41}}$}}%
      \put(6083,2428){\makebox(0,0)[l]{\strut{}${\,h_{22}}$}}%
      \put(6083,812){\makebox(0,0)[l]{\strut{}${\,h_{42}}$}}%
    }%
    \gplgaddtomacro\gplfronttext{%
      \csname LTb\endcsname
      \put(3206,3130){\makebox(0,0){\strut{}$\naZ$}}%
    }%
    \gplbacktext
    \put(0,0){\includegraphics{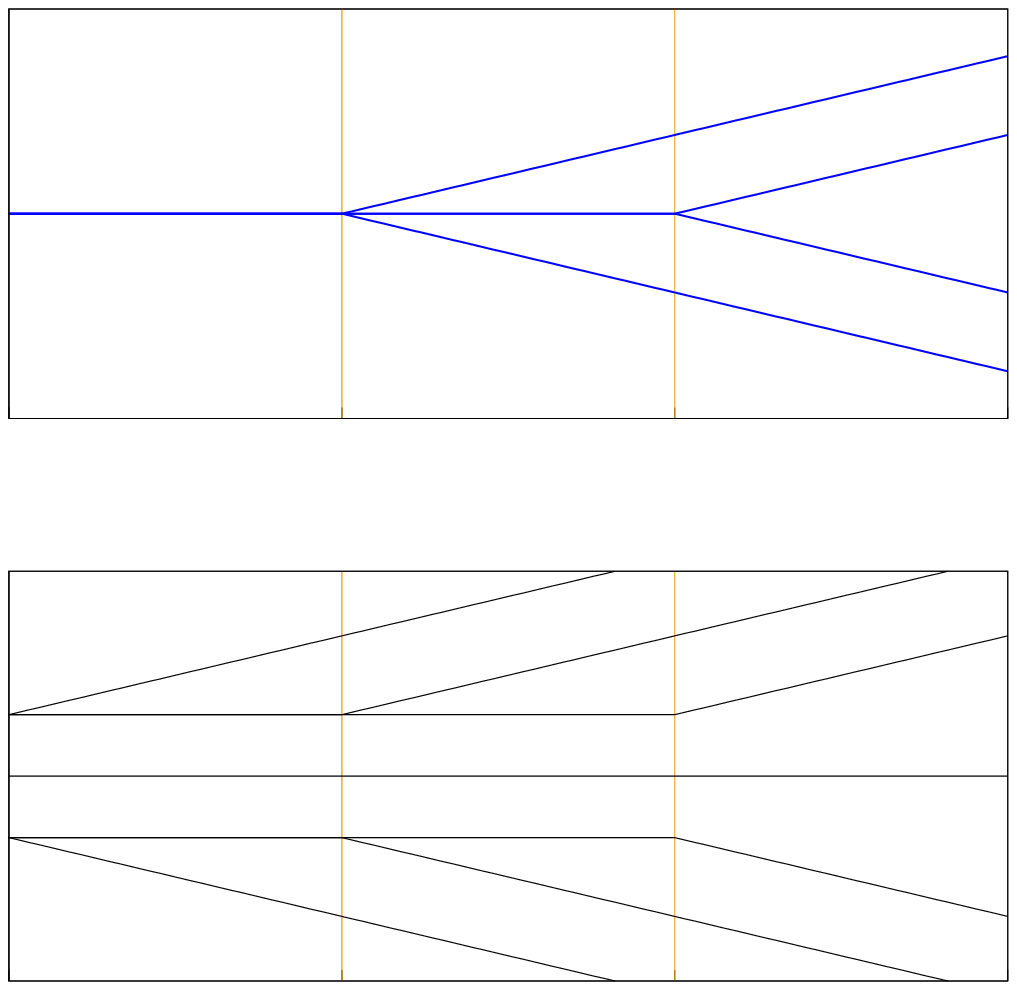}}%
    \gplfronttext
  \end{picture}%
\endgroup

%% file: pic3.tex
\begingroup
  \makeatletter
  \providecommand\color[2][]{%
    \GenericError{(gnuplot) \space\space\space\@spaces}{%
      Package color not loaded in conjunction with
      terminal option `colourtext'%
    }{See the gnuplot documentation for explanation.%
    }{Either use 'blacktext' in gnuplot or load the package
      color.sty in LaTeX.}%
    \renewcommand\color[2][]{}%
  }%
  \providecommand\includegraphics[2][]{%
    \GenericError{(gnuplot) \space\space\space\@spaces}{%
      Package graphicx or graphics not loaded%
    }{See the gnuplot documentation for explanation.%
    }{The gnuplot epslatex terminal needs graphicx.sty or graphics.sty.}%
    \renewcommand\includegraphics[2][]{}%
  }%
  \providecommand\rotatebox[2]{#2}%
  \@ifundefined{ifGPcolor}{%
    \newif\ifGPcolor
    \GPcolorfalse
  }{}%
  \@ifundefined{ifGPblacktext}{%
    \newif\ifGPblacktext
    \GPblacktexttrue
  }{}%
  \let\gplgaddtomacro\g@addto@macro
  \gdef\gplbacktext{}%
  \gdef\gplfronttext{}%
  \makeatother
  \ifGPblacktext
    \def\colorrgb#1{}%
    \def\colorgray#1{}%
  \else
    \ifGPcolor
      \def\colorrgb#1{\color[rgb]{#1}}%
      \def\colorgray#1{\color[gray]{#1}}%
      \expandafter\def\csname LTw\endcsname{\color{white}}%
      \expandafter\def\csname LTb\endcsname{\color{black}}%
      \expandafter\def\csname LTa\endcsname{\color{black}}%
      \expandafter\def\csname LT0\endcsname{\color[rgb]{1,0,0}}%
      \expandafter\def\csname LT1\endcsname{\color[rgb]{0,1,0}}%
      \expandafter\def\csname LT2\endcsname{\color[rgb]{0,0,1}}%
      \expandafter\def\csname LT3\endcsname{\color[rgb]{1,0,1}}%
      \expandafter\def\csname LT4\endcsname{\color[rgb]{0,1,1}}%
      \expandafter\def\csname LT5\endcsname{\color[rgb]{1,1,0}}%
      \expandafter\def\csname LT6\endcsname{\color[rgb]{0,0,0}}%
      \expandafter\def\csname LT7\endcsname{\color[rgb]{1,0.3,0}}%
      \expandafter\def\csname LT8\endcsname{\color[rgb]{0.5,0.5,0.5}}%
    \else
      \def\colorrgb#1{\color{black}}%
      \def\colorgray#1{\color[gray]{#1}}%
      \expandafter\def\csname LTw\endcsname{\color{white}}%
      \expandafter\def\csname LTb\endcsname{\color{black}}%
      \expandafter\def\csname LTa\endcsname{\color{black}}%
      \expandafter\def\csname LT0\endcsname{\color{black}}%
      \expandafter\def\csname LT1\endcsname{\color{black}}%
      \expandafter\def\csname LT2\endcsname{\color{black}}%
      \expandafter\def\csname LT3\endcsname{\color{black}}%
      \expandafter\def\csname LT4\endcsname{\color{black}}%
      \expandafter\def\csname LT5\endcsname{\color{black}}%
      \expandafter\def\csname LT6\endcsname{\color{black}}%
      \expandafter\def\csname LT7\endcsname{\color{black}}%
      \expandafter\def\csname LT8\endcsname{\color{black}}%
    \fi
  \fi
    \setlength{\unitlength}{0.0500bp}%
    \ifx\gptboxheight\undefined%
      \newlength{\gptboxheight}%
      \newlength{\gptboxwidth}%
      \newsavebox{\gptboxtext}%
    \fi%
    \setlength{\fboxrule}{0.5pt}%
    \setlength{\fboxsep}{1pt}%
\begin{picture}(6480.00,6480.00)%
    \gplgaddtomacro\gplbacktext{%
      \csname LTb\endcsname
      \put(330,3460){\makebox(0,0){\strut{}$0$}}%
      \put(3207,3460){\makebox(0,0){\strut{}$1$}}%
      \put(4357,3460){\makebox(0,0){\strut{}$\tau$}}%
      \put(6083,3460){\makebox(0,0){\strut{}$2$}}%
      \put(6083,3680){\makebox(0,0)[l]{\strut{}${\,{\bf v}_1}$}}%
      \put(6083,4860){\makebox(0,0)[l]{\strut{}${\,{\bf v}_2}$}}%
      \put(6083,5251){\makebox(0,0)[l]{\strut{}${\,{\bf v}_3}$}}%
      \put(6083,5449){\makebox(0,0)[l]{\strut{}${\,{\bf v}_4}$}}%
      \put(6083,5567){\makebox(0,0)[l]{\strut{}${\,\bf\ \cdot}$}}%
      \put(6083,5645){\makebox(0,0)[l]{\strut{}${\,\bf\ \cdot}$}}%
      \put(6083,5702){\makebox(0,0)[l]{\strut{}${\,\bf\ \cdot}$}}%
      \put(6083,5744){\makebox(0,0)[l]{\strut{}${\,\bf\ \cdot}$}}%
      \put(6083,5777){\makebox(0,0)[l]{\strut{}${\,\bf\ \cdot}$}}%
      \put(6083,5803){\makebox(0,0)[l]{\strut{}${\,\bf\ \cdot}$}}%
      \put(6083,5824){\makebox(0,0)[l]{\strut{}${\,\bf\ \cdot}$}}%
    }%
    \gplgaddtomacro\gplfronttext{%
      \csname LTb\endcsname
      \put(3206,6369){\makebox(0,0){\strut{}$\VBB$}}%
    }%
    \gplgaddtomacro\gplbacktext{%
      \csname LTb\endcsname
      \put(330,220){\makebox(0,0){\strut{}$0$}}%
      \put(3207,220){\makebox(0,0){\strut{}$1$}}%
      \put(4357,220){\makebox(0,0){\strut{}$\tau$}}%
      \put(6083,220){\makebox(0,0){\strut{}$2$}}%
      \put(6083,440){\makebox(0,0)[l]{\strut{}${\,{\bf u}_1}$}}%
      \put(6083,1620){\makebox(0,0)[l]{\strut{}${\,{\bf u}_2}$}}%
      \put(6083,2012){\makebox(0,0)[l]{\strut{}${\,{\bf u}_3}$}}%
      \put(6083,2210){\makebox(0,0)[l]{\strut{}${\,{\bf u}_4}$}}%
      \put(6083,2328){\makebox(0,0)[l]{\strut{}${\,\bf\ \cdot}$}}%
      \put(6083,2406){\makebox(0,0)[l]{\strut{}${\,\bf\ \cdot}$}}%
      \put(6083,2463){\makebox(0,0)[l]{\strut{}${\,\bf\ \cdot}$}}%
      \put(6083,2505){\makebox(0,0)[l]{\strut{}${\,\bf\ \cdot}$}}%
      \put(6083,2538){\makebox(0,0)[l]{\strut{}${\,\bf\ \cdot}$}}%
      \put(6083,2564){\makebox(0,0)[l]{\strut{}${\,\bf\ \cdot}$}}%
      \put(6083,2585){\makebox(0,0)[l]{\strut{}${\,\bf\ \cdot}$}}%
    }%
    \gplgaddtomacro\gplfronttext{%
      \csname LTb\endcsname
      \put(3206,3130){\makebox(0,0){\strut{}$\UBB$}}%
    }%
    \gplbacktext
    \put(0,0){\includegraphics{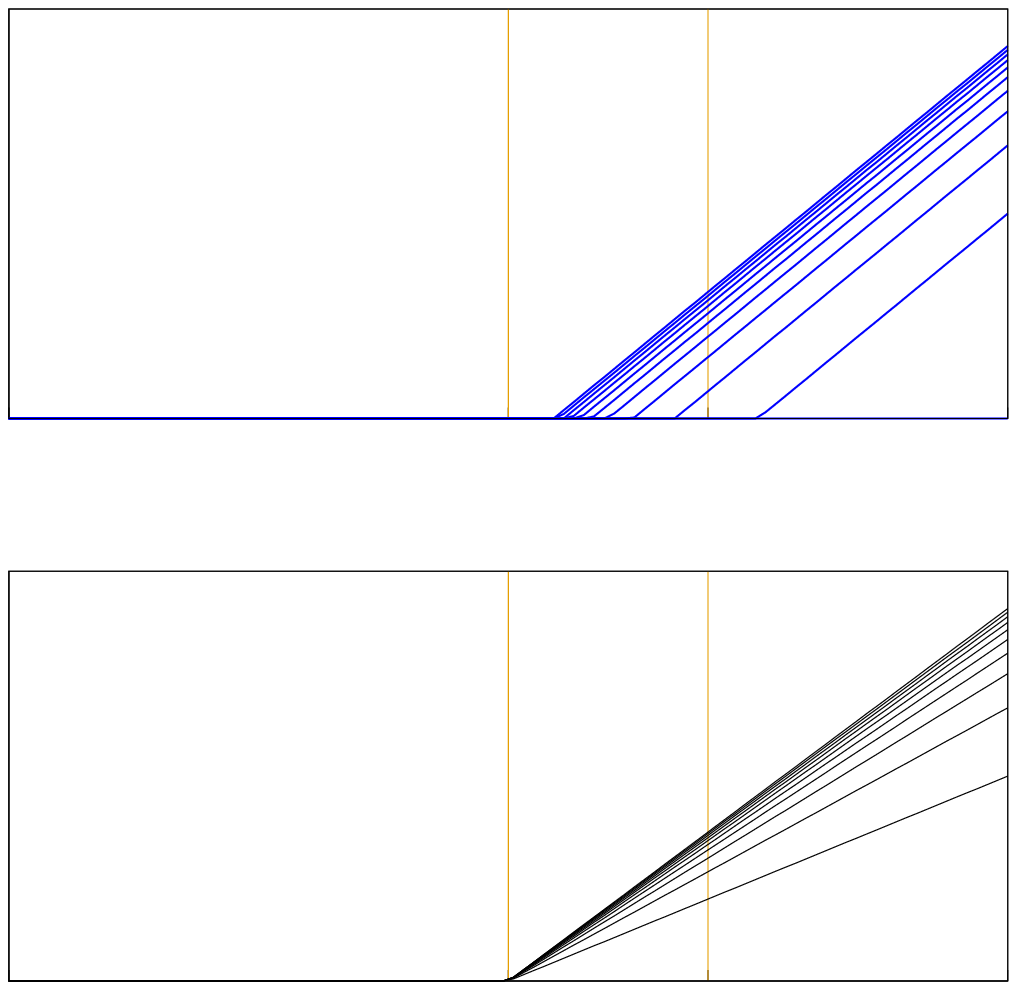}}%
    \gplfronttext
  \end{picture}%
\endgroup

%% file: SerkovDA-DGA-2023.bbl

\begin{thebibliography}{20}
\ifx \bisbn   \undefined \def \bisbn  #1{ISBN #1}\fi
\ifx \binits  \undefined \def \binits#1{#1}\fi
\ifx \bauthor  \undefined \def \bauthor#1{#1}\fi
\ifx \batitle  \undefined \def \batitle#1{#1}\fi
\ifx \bjtitle  \undefined \def \bjtitle#1{#1}\fi
\ifx \bvolume  \undefined \def \bvolume#1{\textbf{#1}}\fi
\ifx \byear  \undefined \def \byear#1{#1}\fi
\ifx \bissue  \undefined \def \bissue#1{#1}\fi
\ifx \bfpage  \undefined \def \bfpage#1{#1}\fi
\ifx \blpage  \undefined \def \blpage #1{#1}\fi
\ifx \burl  \undefined \def \burl#1{\textsf{#1}}\fi
\ifx \doiurl  \undefined \def \doiurl#1{\url{https://doi.org/#1}}\fi
\ifx \betal  \undefined \def \betal{\textit{et al.}}\fi
\ifx \binstitute  \undefined \def \binstitute#1{#1}\fi
\ifx \binstitutionaled  \undefined \def \binstitutionaled#1{#1}\fi
\ifx \bctitle  \undefined \def \bctitle#1{#1}\fi
\ifx \beditor  \undefined \def \beditor#1{#1}\fi
\ifx \bpublisher  \undefined \def \bpublisher#1{#1}\fi
\ifx \bbtitle  \undefined \def \bbtitle#1{#1}\fi
\ifx \bedition  \undefined \def \bedition#1{#1}\fi
\ifx \bseriesno  \undefined \def \bseriesno#1{#1}\fi
\ifx \blocation  \undefined \def \blocation#1{#1}\fi
\ifx \bsertitle  \undefined \def \bsertitle#1{#1}\fi
\ifx \bsnm \undefined \def \bsnm#1{#1}\fi
\ifx \bsuffix \undefined \def \bsuffix#1{#1}\fi
\ifx \bparticle \undefined \def \bparticle#1{#1}\fi
\ifx \barticle \undefined \def \barticle#1{#1}\fi
\bibcommenthead
\ifx \bconfdate \undefined \def \bconfdate #1{#1}\fi
\ifx \botherref \undefined \def \botherref #1{#1}\fi
\ifx \url \undefined \def \url#1{\textsf{#1}}\fi
\ifx \bchapter \undefined \def \bchapter#1{#1}\fi
\ifx \bbook \undefined \def \bbook#1{#1}\fi
\ifx \bcomment \undefined \def \bcomment#1{#1}\fi
\ifx \oauthor \undefined \def \oauthor#1{#1}\fi
\ifx \citeauthoryear \undefined \def \citeauthoryear#1{#1}\fi
\ifx \endbibitem  \undefined \def \endbibitem {}\fi
\ifx \bconflocation  \undefined \def \bconflocation#1{#1}\fi
\ifx \arxivurl  \undefined \def \arxivurl#1{\textsf{#1}}\fi
\csname PreBibitemsHook\endcsname

\bibitem[\protect\citeauthoryear{Fleming}{1961}]{Fleming:1961}
\begin{barticle}
\bauthor{\bsnm{Fleming}, \binits{W.H.}}:
\batitle{{The Convergence Problem for Differential Games}}.
\bjtitle{Journal of Mathematical Analysis and Applications}
\bvolume{3}(\bissue{1}),
\bfpage{102}--\blpage{116}
(\byear{1961})
\doiurl{10.1016/0022-247X(61)90009-9}
\end{barticle}
\endbibitem

\bibitem[\protect\citeauthoryear{Pontryagin}{1967}]{Pontr1967et}
\begin{barticle}
\bauthor{\bsnm{Pontryagin}, \binits{L.S.}}:
\batitle{Linear differential games.~2}.
\bjtitle{Soviet Math. Dokl.}
\bvolume{8},
\bfpage{910}--\blpage{912}
(\byear{1967})
\end{barticle}
\endbibitem

\bibitem[\protect\citeauthoryear{Blagodatskikh and
  Petrov}{2019}]{BlaPet-DGA2019}
\begin{barticle}
\bauthor{\bsnm{Blagodatskikh}, \binits{A.I.}},
\bauthor{\bsnm{Petrov}, \binits{N.N.}}:
\batitle{Simultaneous multiple capture of rigidly coordinated evaders}.
\bjtitle{Dyn Games Appl}
\bvolume{9},
\bfpage{594}--\blpage{613}
(\byear{2019})
\doiurl{10.1007/s13235-019-00300-8}
\end{barticle}
\endbibitem

\bibitem[\protect\citeauthoryear{Chernov}{2014}]{Chernov2014ARC}
\begin{barticle}
\bauthor{\bsnm{Chernov}, \binits{A.V.}}:
\batitle{On {Volterra} functional operator games on a given set}.
\bjtitle{Autom Remote Control}
\bvolume{75},
\bfpage{787}--\blpage{803}
(\byear{2014})
\doiurl{10.1134/S0005117914040195}
\end{barticle}
\endbibitem

\bibitem[\protect\citeauthoryear{Petrosyan and Zenkevich}{1996}]{PetZen1996et}
\begin{bbook}
\bauthor{\bsnm{Petrosyan}, \binits{L.A.}},
\bauthor{\bsnm{Zenkevich}, \binits{N.A.}}:
\bbtitle{Game Theory. {Transl}. from the {Russian} by {J}. {M}. {Donetz}}.
\bpublisher{Singapore: World Scientific}, 
(\byear{1996})
\end{bbook}
\endbibitem

\bibitem[\protect\citeauthoryear{Khlopin and Chentsov}{2005}]{CheKhlo2005DU}
\begin{barticle}
\bauthor{\bsnm{Khlopin}, \binits{D.V.}},
\bauthor{\bsnm{Chentsov}, \binits{A.G.}}:
\batitle{{On a Control Problem with Incomplete Information: Quasistrategies and
  Control Procedures with a Model}}.
\bjtitle{Differ.Uravn.}
\bvolume{41}(\bissue{12}),
\bfpage{1727}--\blpage{1742}
(\byear{2005})
\end{barticle}
\endbibitem

\bibitem[\protect\citeauthoryear{Gomoyunov and Serkov}{2021}]{GomSer-UDSU2021}
\begin{barticle}
\bauthor{\bsnm{Gomoyunov}, \binits{M.I.}},
\bauthor{\bsnm{Serkov}, \binits{D.A.}}:
\batitle{On guarantee optimization in control problem with finite set of
  disturbances}.
\bjtitle{Vestnik Udmurtskogo Universiteta. Matematika. Mekhanika. Komp'yuternye
  Nauki}
\bvolume{31}(\bissue{4}),
\bfpage{613}--\blpage{628}
(\byear{2021})
\doiurl{10.35634/vm210406}
\end{barticle}
\endbibitem

\bibitem[\protect\citeauthoryear{Chentsov}{2001a}]{Chentsov2001DEI}
\begin{barticle}
\bauthor{\bsnm{Chentsov}, \binits{A.G.}}:
\batitle{{Nonanticipating Multimappings and Their Construction by the Method of
  Program Iterations: I.}}
\bjtitle{Differential Equations}
\bvolume{37},
\bfpage{498}--\blpage{509}
(\byear{2001})
\doiurl{10.1023/A:1019275422741}
\end{barticle}
\endbibitem

\bibitem[\protect\citeauthoryear{Chentsov}{2001b}]{Chentsov2001DEII}
\begin{barticle}
\bauthor{\bsnm{Chentsov}, \binits{A.G.}}:
\batitle{{Nonanticipating Multimappings and Their Construction by the Method of
  Program Iterations: II}}.
\bjtitle{Differential Equations}
\bvolume{37},
\bfpage{713}--\blpage{723}
(\byear{2001})
\doiurl{10.1023/A:1019224800877}
\end{barticle}
\endbibitem

\bibitem[\protect\citeauthoryear{Serkov}{2017}]{Ser_UDSU2017}
\begin{barticle}
\bauthor{\bsnm{Serkov}, \binits{D.A.}}:
\batitle{Unlocking of predicate: application to constructing a non-anticipating
  selection}.
\bjtitle{Vestnik Udmurtskogo Universiteta. Matematika. Mekhanika. Komp'yuternye
  Nauki}
\bvolume{27}(\bissue{2}),
\bfpage{283}--\blpage{291}
(\byear{2017})
\doiurl{10.20537/vm170211}
\end{barticle}
\endbibitem

\bibitem[\protect\citeauthoryear{Chentsov}{1978}]{Chentsov1978dep-e}
\begin{botherref}
\oauthor{\bsnm{Chentsov}, \binits{A.G.}}:
Selections of multivalued strategies in differential games.
Manuscript deposited at VINITI 3101-78Dep,
IMM USC of the USSR Academy of Sciences, Sverdlovsk,
Sverdlovsk
(1978).
(in Russian)
\end{botherref}
\endbibitem

\bibitem[\protect\citeauthoryear{Cardaliaguet and
  Plaskacz}{2000}]{CARDALIAGUET-PLASKACZ-SIAM2000}
\begin{barticle}
\bauthor{\bsnm{Cardaliaguet}, \binits{P.}},
\bauthor{\bsnm{Plaskacz}, \binits{S.}}:
\batitle{Invariant solutions of differential games and
  {Hamilton--Jacobi--Isaacs} equations for time-measurable hamiltonians}.
\bjtitle{SIAM J. Control Optim}
\bvolume{38}(\bissue{5}),
\bfpage{1501}--\blpage{1520}
(\byear{2000})
\doiurl{10.1137/S0363012998296219}
\end{barticle}
\endbibitem

\bibitem[\protect\citeauthoryear{Serkov and
  Chentsov}{2020}]{CheSer_TRIMM2019et}
\begin{barticle}
\bauthor{\bsnm{Serkov}, \binits{D.A.}},
\bauthor{\bsnm{Chentsov}, \binits{A.G.}}:
\batitle{{On the Construction of a Nonanticipating Selection of a Multivalued
  Mapping}}.
\bjtitle{Proceedings of the Steklov Institute of Mathematics}
\bvolume{309}(\bissue{Suppl.1}),
\bfpage{125}--\blpage{138}
(\byear{2020})
\doiurl{10.1134/S008154382004015X}
\end{barticle}
\endbibitem

\bibitem[\protect\citeauthoryear{Chentsov}{1997}]{Chentsov1997DAN}
\begin{barticle}
\bauthor{\bsnm{Chentsov}, \binits{A.G.}}:
\batitle{The iterative realization of nonanticipating multivalued mappings}.
\bjtitle{Doklady Mathematics}
\bvolume{56}(\bissue{3}),
\bfpage{927}--\blpage{930}
(\byear{1997})
\end{barticle}
\endbibitem

\bibitem[\protect\citeauthoryear{Subbotin and Chentsov}{1981}]{SubChe81e}
\begin{bbook}
\bauthor{\bsnm{Subbotin}, \binits{A.I.}},
\bauthor{\bsnm{Chentsov}, \binits{A.G.}}:
\bbtitle{Optimization of Guarantee in Control Problems},
p. \bfpage{288}.
\bpublisher{Nauka},
\blocation{M.}
(\byear{1981}).
\bcomment{(in Russian)}
\end{bbook}
\endbibitem

\bibitem[\protect\citeauthoryear{Serkov}{2022}]{Ser_CTMM2022e}
\begin{bchapter}
\bauthor{\bsnm{Serkov}, \binits{D.A.}}:
\bctitle{Step-by-step construction of optimal motion and non-anticipative
  multi-selectors}.
In: \bbtitle{Control Theory and Mathematical Modeling: Vseros. Conf. with
  International Participation, Dedicated Memory of Prof. N.V. Azbelev and Prof.
  E.L. Tonkov (STMM 2022), Izhevsk, June 13-17, 2022 : Materials},
pp. \bfpage{219}--\blpage{223}
(\byear{2022})
\end{bchapter}
\endbibitem

\bibitem[\protect\citeauthoryear{Engelking}{1985}]{Engelking1986e}
\begin{bbook}
\bauthor{\bsnm{Engelking}, \binits{R.}}:
\bbtitle{General Topology},
p. \bfpage{752}.
\bpublisher{Panstwowe Wydawnictwo Naukowe},
\blocation{Warszawa}
(\byear{1985})
\end{bbook}
\endbibitem

\bibitem[\protect\citeauthoryear{Kuratowski}{1966}]{Kura1966e}
\begin{bbook}
\bauthor{\bsnm{Kuratowski}, \binits{K.}}:
\bbtitle{Topology. {Volume I}}.
\bpublisher{Academic Press},
\blocation{NY and London}
(\byear{1966}).
\doiurl{10.1016/C2013-0-11022-7}
\end{bbook}
\endbibitem

\bibitem[\protect\citeauthoryear{Krasovskij et~al.}{1972}]{NNK:143et}
\begin{barticle}
\bauthor{\bsnm{Krasovskij}, \binits{N.N.}},
\bauthor{\bsnm{Subbotin}, \binits{A.I.}},
\bauthor{\bsnm{Ushakov}, \binits{V.N.}}:
\batitle{A minimax differential game}.
\bjtitle{Sov. Math., Dokl.}
\bvolume{13},
\bfpage{1200}--\blpage{1204}
(\byear{1972})
\end{barticle}
\endbibitem

\bibitem[\protect\citeauthoryear{Pshenichnyj}{1969}]{Pshenichny:1969et}
\begin{barticle}
\bauthor{\bsnm{Pshenichnyj}, \binits{B.N.}}:
\batitle{The structure of differential games}.
\bjtitle{Sov. Math., Dokl.}
\bvolume{10},
\bfpage{70}--\blpage{72}
(\byear{1969})
\end{barticle}
\endbibitem

\end{thebibliography}
